\documentclass{amsproc}

\usepackage[swedish,english]{babel}
\usepackage[utf8]{inputenc}
\usepackage[T1]{fontenc}
\usepackage{amsmath,amssymb,amsfonts}
\usepackage{url}
\usepackage[dvipsnames,usenames]{color}
\usepackage{array}
\usepackage{accents}
\usepackage{tabu}
\usepackage[normalem]{ulem}
\usepackage{tikz}

\newcolumntype{L}{>{\displaystyle}l}
\newcolumntype{C}{>{\displaystyle}c}
\newcolumntype{R}{>{\displaystyle}r}

\renewcommand{\tfrac}{\genfrac{}{}{}1}

\newcommand{\R}{\mathbb R}
\newcommand{\N}{\mathbb N}
\newcommand{\Z}{\mathbb Z}
\newcommand{\C}{\mathbb C}

\newcommand{\cont}{\mathcal C}
\newcommand{\D}{\mathbb D}
\newcommand{\T}{\mathbb T}
\renewcommand{\Im}{\mathrm{Im}}
\renewcommand{\Re}{\mathrm{Re}}
\def\E{{\mathrm{e}}}

\def\til{\widetilde}

\newcommand{\pois}{\mathcal{P}}
\newcommand{\four}{\mathcal{F}}
\newcommand{\schw}{\mathcal{S}}

\def\I{\mathfrak{i}}
\newcommand{\diff}{\mathrm{d}}
\renewcommand{\bar}{\overline}

\newcommand{\supp}{\mathrm{supp}}
\newcommand{\sgn}{\mathrm{sgn}}
\newcommand{\HN}{Herglotz-Nevan\-linna }
\newcommand{\ntto}{\:\scriptsize{\xrightarrow{\vee\:}}\:}
\newcommand{\Leb}{\mathrm{L}}
\newcommand{\Log}{\mathrm{Log}}

\newcommand{\ie}{\textit{i.e.}\/ }
\newcommand{\eg}{\textit{e.g.}\/ }
\newcommand{\cf}{\textit{cf.}\/ }

\usepackage{bm}
\renewcommand{\vec}{\bm}
\newcommand{\mat}[1]{#1}
\newcommand{\tilvec}[1]{\til{\vec{#1}}}

\allowdisplaybreaks

\theoremstyle{definition} 
\newtheorem{define}{Definition}[section]
\newtheorem{example}[define]{Example}
\newtheorem{remark}[define]{Remark}

\theoremstyle{plain} 
\newtheorem{lemma}[define]{Lemma}
\newtheorem{thm}[define]{Theorem}
\newtheorem{prop}[define]{Proposition}
\newtheorem{coro}[define]{Corollary}

\newcommand{\NL}{{\mathbf{M}_{{\rm sep}}}}
\newcommand{\NLP}{{\mathbf{M}_{{\rm sep}}^+}}

\numberwithin{equation}{section}

\begin{document}

\title[On the structure of Nevanlinna measures]{On the structure of Nevanlinna measures}

\author{Mitja Nedic}
\address{Mitja Nedic, Department of Mathematics and Statistics, University of Helsinki, PO Box 68, FI-00014 Helsinki, Finland, orc-id: 0000-0001-7867-5874}
\curraddr{}
\email{mitja.nedic@helsniki.fi}

\author{Eero Saksman}
\address{Eero Saksman, Department of Mathematics and Statistics, University of Helsinki, PO Box 68, FI-00014 Helsinki, Finland.}
\curraddr{}
\email{eero.saksman@helsniki.fi}


\begin{abstract}
In this paper, we study the structural properties of Nevanlinna measures, \ie Borel measures that arise in the integral representation of Herglotz-Nevanlinna functions. In particular, we give a characterization of these measures in terms of their Fourier transform, characterize measures supported on hyperplanes including extremal measures, describe the structure of the singular part of the measures when some variable are set to a fixed value, and provide estimates for the measure of expanding and shrinking cubes. Corresponding results are stated also in the setting of the polydisc where applicable, and some of our proofs are actually perfomed via the polydisc.
\end{abstract}

\maketitle

\section{Introduction}

When considering holomorphic functions in one or several variables, those that map a given domain into a half-plane play have their on special interest, and they play a special role in many areas and applications, both within and outside mathematics. For functions of one variable, a few examples of such applications are found within extension theory of symmetric operators \cite{AkhiezerGlazman1993,LangerTextorius1977} or spectral theory of Sturm-Liouville problems and perturbations \cite{Aronszajn1957,AronszajnBrown1970,Donoghue1965,KacKrein1974}, and when describing representations or realizations of passive and non-passive systems \cite{BernlandEtal2011,IvanenkoETAL2019,IvanenkoETAL2020,Zemanian1963,Zemanian1965}. For applications of such functions in several variables, we mention the theory of operator monotone functions \cite{AglerEtal2012} or representations of multidimensional passive systems \cite{Vladimirov1979}. 

With regards to the domain of the function in question, one usually restricts themselves to one of the following the two cases. Firstly, in the poly-upper half-plane, where we consider the class of \HN functions, \ie functions with non-negative imaginary part, see Section \ref{subsec:HN_functions}. Secondly, in the unit polydisc, where we consider RP-functions, \ie functions with non-negative real part, see Section \ref{subsec:RP_functions}. We are following here the tradition of the literature, although one may  note that functions on the unit polydisc with non-negative imaginary part are considered instead in \eg \cite{VladimirovDrozzinov1974}, and more general domains are considered in \eg \cite{AizenbergDautov1976,AizenbergYuzhakov1983,Vladimirov1969}.

Both \HN functions and RP-functions can be characterized via an integral representations formula, see Theorems \ref{thm:intRep_Nvar} and \ref{thm:RP_intRep_Nvar}, respectively. In short, both representations are of the form 
$$L(\vec{z}) + \int_D \mathcal{K}_n(\vec{z},\vec{t}) \diff\rho(\vec{t}),$$
where, $L$ denotes a linear term and $\mathcal{K}_n$ is a kernel function, both depending on the domain in question, while $D$ is the distinguished boundary of the domain and $\rho$ is a positive Borel measure on $D$ satisfying certain conditions. It is these classes of measures, called \emph{Nevanlinna measures} in the case of the poly-upper half-plane, \cf Section \ref{subsec:HN_functions}, or \emph{RP-measures} in the case of the polydisc, \cf Section \ref{subsec:RP_functions}, that are of particular interest as they encode many essential properties of the function they represent via the respective integral representation formula.

In particular, the goal of this paper is to describe different structural aspects of Nevanlinna measures. First, we focus on the Fourier transform of Nevanlina measures and prove, in Theorem \ref{thm:Fourier_char}, that a measure $\mu$ is a Nevanlinna measure if and only if its Fourier transform, understood as a Schwartz distribution, satisfies
$$\supp(\widehat\mu) \; \subset\; [0,\infty)^n\cup (\infty,0]^n.$$
We will give several applications of this result. Especially, we describe all measures supported on certain hyperplanes in $\R^n$ while at the same time being able to discern which of these measures are extremal, \cf Theorem \ref{thm:hyperplanes} and Corollary \ref{coro:measure_on_lines}. Furthermore, in Theorem \ref{thm:lines_with_positive_slope}, we show that the restriction of an arbitrary Nevanlinna measure to a certain hyperplane is always zero.

Second, we  study how the singular part of a Nevanlinna measure behaves when some of the variables of a function are fixed. Our main result in this direction is presented in Theorem \ref{thm:singular_part_HN}, which states that if the measure has a non-trivial singular part at one fixed value, then it has one at all values. Moreover, this singular part turns out to be independent of the fixed variable, allowing for a particular decomposition of the function in question. Theorem \ref{thm:singular_part_polydisc} presents the same result in the language of RP-functions and RP-measures. Finally, we investigate estimates on the measure of expanding cubes in Propositions \ref{prop:volume_estimate} and \ref{prop:limsup}.

The structure of the paper is as follows. After the introduction we focus in Section \ref{sec:prerec} on giving a detailed overview of the prerequisites regarding \HN functions and RP-functions as well as their representing measure. Section \ref{sec:Fourier} is devoted to the main result regarding the characterization of Nevanlinna measures via the Fourier transform as well as several corollaries and examples. Section \ref{sec:fixed_var} concentrates on the problem where some of the independent variables are fixed, presenting the main results on both the polydisc and poly-upper half-plane. Finally, Section \ref{sec:estimates} contains some results and examples regarding estimates on the measures of cubes. In each section we provide references to  previous works on related topics.

\section{Holomorphic functions having non-negative imaginary or real part}
\label{sec:prerec}

In this section, we first review the necessary prerequisites regarding \HN functions and Nevanlinna measures, followed by RP-functions and RP-measures. Finally, we describe in detail how the two classes and their respective integral representations relate to each other.

\subsection{\HN functions in the poly-upper half-plane}\label{subsec:HN_functions}

In the poly-upper half-plane $\C^{+n}$, \ie 
$$\C^{+n}:= (\C^+)^n = \big\{\vec{z}\in\C^n \,\big |\,\forall j=1,2,\ldots, n:  \Im[z_j]>0 \big\},$$
we consider the following class of functions, \cf \cite{LugerNedic2017,LugerNedic2019,Vladimirov1969,Vladimirov1979}.

\begin{define}\label{def:HN_Nvar}
A function $q\colon \C^{+n} \to \C$ is called a \emph{\HN function} if it is holomorphic with a non-negative imaginary part.
\end{define}

We recall also the integral representation theorem for \HN functions of several variables, \cf \cite[Thm. 4.1 and Thm. 5.1]{LugerNedic2019}.

\begin{thm}\label{thm:intRep_Nvar}
A function $q\colon \C^{+n} \to \C$ is a \HN function if and only if $q$ can be written, for every $\vec{z} \in \C^{+n}$, as
\begin{equation}\label{eq:intRep_Nvar}
q(\vec{z}) = a + \sum_{\ell=1}^nb_\ell z_\ell + \frac{1}{\pi^n}\int_{\R^n}K_n(\vec{z},\vec{t})\diff\mu(\vec{t}),
\end{equation}
where $a \in \R$, $\vec{b} \in [0,\infty)^n$, the kernel $K_n\colon\C^{+n} \times \R^n \to \C$ is defined as
$$K_n(\vec{z},\vec{t}) := \I\left(\frac{2}{(2\I)^n}\prod_{\ell=1}^n\left(\frac{1}{t_\ell-z_\ell}-\frac{1}{t_\ell+\I}\right)-\frac{1}{(2\I)^n}\prod_{\ell=1}^n\left(\frac{1}{t_\ell-\I}-\frac{1}{t_\ell+\I}\right)\right)$$
and $\mu$ is a positive Borel measure on $\R^n$ satisfying the growth condition
\begin{equation}
\label{eq:measure_growth}
\int_{\R^n}\prod_{\ell=1}^n\frac{1}{1+t_\ell^2}\diff\mu(\vec{t}) < \infty
\end{equation}
and the Nevanlinna condition, \ie
\begin{equation}
    \label{eq:measure_Nevan}
    \int_{\R^n}\frac{1}{(t_{j_1}-z_{j_1})^2(t_{j_2}-\bar{z_{j_2}})^2}\:\prod_{\substack{\ell=1 \\ \ell \neq j_1,j_2}}^n\left(\frac{1}{t_\ell-z_\ell} - \frac{1}{t_\ell-\bar{z_\ell}}\right)\diff\mu(\vec{t}) = 0
\end{equation}
for all $\vec{z} \in \C^{+n}$ and for all indices $j_1,j_2 \in \N$, such that $1 \leq j_1 < j_2 \leq n$. Furthermore, for a given function $q$, the triple of representing parameters $(a,\vec{b},\mu)$ is unique.
\end{thm}

A positive Borel measure $\mu$ on $\R^n$ satisfying conditions \eqref{eq:measure_growth} and \eqref{eq:measure_Nevan} is called \emph{ a Nevanlinna measure}, see \eg \cite{LugerNedic2019,LugerNedic2021,Nedic2020}. By \cite[Cor. 4.6(ii)]{LugerNedic2019}, it holds for any Nevanlinna measure $\mu$ that
$$\int_{\R^n}\Im[K_n(\vec{z},\vec{t})]\diff\mu(\vec{t}) = \int_{\R^n}\pois_n(\vec{z},\vec{t})\diff\mu(\vec{t}),$$
where $\pois_n$ denotes the Poisson kernel of the poly-upper half-plane, \ie
\begin{equation}
    \label{eq:Poisson_half_plane}
    \pois_n(\vec{z},\vec{t}) := \prod_{j=1}^n\frac{\Im[z_j]}{|z_j - t_j|^2}.
\end{equation}
Equivalently, Nevanlinna measures may be described as precisely those positive Borel measures on $\R^n$ satisfying condition \eqref{eq:measure_growth} for which the function
$$\vec{z} \mapsto \int_{\R^n}\pois_n(\vec{z},\vec{t})\diff\mu(\vec{t})$$
is pluriharmonic on $\C^{+n}$.

An important consequence of Theorem \ref{thm:intRep_Nvar} says that the measure $\mu$ may be recovered form the function $q$ via the Stieltjes inversion formula, \cf \cite[Prop. 4.1]{LugerNedic2017} and \cite[Cor. 4.6(viii)]{LugerNedic2019}.

\begin{prop}
Let $q$ be a \HN function and $\mu$ its representing measure in the sense of Theorem \ref{thm:intRep_Nvar}. Let $\psi\colon \R^n \to \R$ be a $\cont^1$-function for which there exists a constant $C \geq 0$ such that
$|\psi(\vec{x})| \leq C\,\prod_{j=1}^n(1+x_j^2)^{-1}$ for all $\vec{x} \in \R^n$. Then, it holds that
\begin{equation}\label{eq:Stieltjes_inversion}
\int_{\R^n}\psi(\vec{t})\diff\mu(\vec{t})=\lim\limits_{\vec{y} \to \vec{0}^+}\int_{\R^n}\psi(\vec{x})\Im[q(\vec{x} + \I\:\vec{y})]\diff \vec{x}.
\end{equation}
\end{prop}

In other words, the measure $\mu$ can be thought of as the limit of the function $\vec{z} \mapsto \Im[q(z)]$ when $\vec{z}$ non-tangentially approaches the distinguished boundary of the poly-upper half-plane. This limit need not exist at every point of the distinguished boundary, but exists in a distributional sense as described in formula \eqref{eq:Stieltjes_inversion}. For the readers benefit we recall couple of interesting examples of Nevanlinna measures in $\lambda_{\R^2}$.  
The functions
$$(z_1,z_2) \mapsto \frac{-1}{z_1+z_2+\I} \quad\text{and}\quad (z_1,z_2) \mapsto \sqrt{z_1}\sqrt{z_2},$$
where in the second example the square root is taken to have a branch cut along the negative part of the real axis,  lead to absolutely continuous Nevanlinna mesures, as is easily verified. In turn, $-(z_1+z_2)^{-1}$ has singular Nevanlinna measure, supported on line $\{x_1+x_2=0\}$.

Finally, we would like to mention the following properties satisfied by a Nevanlinna measure $\mu$ on $\R^n$ with $n \geq 2$.
\begin{itemize}
    \item{The measure $\mu$ cannot be finite unless it is identically zero, \cf \cite[Prop. 4.3]{LugerNedic2017}, see also Corollary \ref{coro:finite_measure} for a proof of this statement using the results of the present paper. When $n = 1$, there is no such restriction.}
    \item{The restriction of $\mu$ to any hyperplane that is orthogonal against one of the coordinate axes is a multiple of the Lebesgue measure in dimension $n-1$, \cf \cite[Thm. 3.4]{LugerNedic2021}  (see also Lemma \ref{lem:hyperplanes} in the case of the polydisc below). In particular, all points have zero mass, \cf \cite[Prop. 4.4]{LugerNedic2017} or \cite[Cor. 3.8]{LugerNedic2021}. Theorems \ref{thm:hyperplanes} and \ref{thm:lines_with_positive_slope} below complete the picture by characterizing all Nevanlinna measures supported on hyperplanes not parallel to a coordinate axis. When $n=1$, there are no hyperplanes to consider and point masses are possible.}
    \item{The (topological) support of $\mu$ obeys certain geometric restrictions, \cf \cite[Thms. 3.11, 3.17 and 3.25]{LugerNedic2021}. When $n=1$, any Borel subset of $\R$ can appear as the support of some Nevanlinna measure.}
    \item{If $\mu$ can be written as a product measure, then at least one of the factors must be equal to a constant multiple of the Lebesgue measure, \cf \cite[Thm. 4.1]{Nedic2020}, see also Corollary \ref{coro:product_measure} for a proof of this statement using the results of the present paper. When $n=1$, there are no product measures to consider.}
\end{itemize}

\subsection{RP-functions in the unit polydisc}\label{subsec:RP_functions}

In the unit polydisc $\D^n$, the following class of functions plays an analogous role to \HN functions in the poly-upper half-plane, \cf \cite{KoranyiPukanszky1963,McDonald1987,Rudin1969,Zemanian1965}.

\begin{define}\label{def:RP_functions}
A function $G\colon \D^n \to \C$ is called a \emph{RP-function} if it is holomorphic with a non-negative real part.
\end{define}

An analogous version of Theorem \ref{thm:intRep_Nvar} for RP-functions is the following \cite[Thm. 1]{KoranyiPukanszky1963}.

\begin{thm}\label{thm:RP_intRep_Nvar}
A function $G\colon \D^n \to \C$ is a RP-function if and only if $G$ can be written, for every $\vec{z} \in \D^n$, as
\begin{equation}\label{eq:RP_intRep_Nvar}
G(\vec{z}) = \I\,A + \int_{\T^{n}} \left(-1+2\prod_{\ell=1}^{n}\frac{1}{1-z_\ell\,\overline{w_\ell}}\right) \diff\nu(\vec{w}),
\end{equation}
where $A \in \R$ and $\nu$ is a finite positive Borel measure on $\T^n$ such that
\begin{equation}\label{eq:RP_condition}
    \int_{\T^n}w_1^{j_1}\ldots w_n^{j_n}\diff\nu(\vec{w}) = 0
\end{equation}
for every multi-index $\vec{j} \in \Z^n$ with at least one positive and one negative entry. Moreover, it holds that $A = \Im[G(\vec{0})]$ and that for every RP-function $G$ there exists a unique measure $\nu$ such that formula \eqref{eq:RP_intRep_Nvar} holds for all $\vec{z} \in \D^n$.
\end{thm}

In the literature, there exist different established ways of writing representation \eqref{eq:RP_intRep_Nvar} based on how the torus $\T^n$ and the measure $\nu$ are parametrized. In particular, we wish to recall the following variant. By parametrizing the torus as $[0,2\pi)^n$ and applying the change of variables $\,w = \E^{\I\,s}$ between $w \in \T$ and $s \in [0,2\pi)$ in each variable separately, one may re-parametrize the measure $\til{\nu}$ as a finite positive Borel measure $\sigma$ on $[0,2\pi)^n$. In particular, using this notation, it holds for any Borel function $f \in \Leb^1(\T^n,\nu)$ that
$$\int_{\T^n}f(\vec{w})\diff\nu(\vec{w}) =  \frac{1}{(2\pi)^n}\int_{[0,2\pi)^n}f(\E^{\I\,s_1},\ldots,\E^{\I\,s_n})\diff\sigma(\vec{s}).$$
For convenience, we write
\begin{equation}
    \label{eq:nu_to_sigma}
    \sigma = \Sigma(\nu) \quad\text{or}\quad \nu = \Sigma^{-1}(\sigma)
\end{equation}
whenever we are referring to two measures related by exactly the reparametrizations reviewed above.

We introduce now \emph{RP-measures} as finite positive Borel measures on $\T^n$ satisfying condition \eqref{eq:RP_condition}, see \eg \cite{Ahern1973,LugerNedic2017,LugerNedic2021,McDonald1982,McDonald1986,McDonald1990,Nedic2020,Rudin1969}. We remark that these measures are also known under the name \emph{measures with vanishing mixed Fourier coefficients} as condition \eqref{eq:RP_condition} may be written as $\hat{\nu}(\vec{j}) = 0$ for $\vec{j} \in \Z^n$ as before, where, as common, $\hat{\nu}(\vec{j})$ denotes the $\vec{j}$-th Fourier coefficient of $\nu$, \ie
$$\hat{\nu}(\vec{j}) := \int_{\T^n}w_1^{-j_1}\ldots w_n^{-j_n}\diff\nu(\vec{w}).$$
RP-measures share an intricate connection to the Poisson kernel of $\D^n$ in analogous way as Nevanlinna measures are connected to the Poisson kernel of $\C^{+n}$. In particular, it holds for any RP-function $G$ with representing measure $\nu$ that
$$\Re[G(\vec{z})] = \int_{\T^n}\left(-1+2\,\Re\left[\prod_{\ell=1}^{n}\frac{1}{1-z_\ell\overline{w_\ell}}\right]\right)\diff\nu(\vec{w}) = \int_{\T^n}P_n(\vec{z},\vec{w})\diff\nu(\vec{w}),$$
where $\vec{z} \in \D^n$ and $P_n$ denotes the Poisson kernel of the unit polydisc, \ie
$$P_n(\vec{z},\vec{w}) := \prod_{\ell=1}^n\Re\left[\frac{w_\ell + z_\ell}{w_\ell - z_\ell}\right] \;=\; \prod_{\ell=1}^n \frac{1-|z_\ell|^2}{|w_\ell-z_\ell|^2}.$$

We note also that the set of RP-measures constitutes a convex cone within the Banach space of all complex Borel measures on $\T^n$, where the norm is given by the total variation. Hence, by the Krein-Milman theorem \cite[V.8.4]{DunfordSchwartz1958}, the set of of all RP-measures of total variation $\leq 1$ is equal to the convex hull of its extremal points, \ie points that cannot be written as a non-trivial convex combination of two points form the same set. In dimension 1, since every finite Borel measure is a RP-measure, these extremal measures are precisely Dirac measures supported at the different points of $\T^1$, while in higher dimensions, no complete description currently exists.

In analogy with Nevanlinna measures on $\R^n$, the following properties are satisfied by a RP-measure $\nu$ on $\T^n$ with $n \geq 2$.

\begin{itemize}
    \item{The restriction of $\nu$ to any hyperplane passing through a given point is equal to a multiple of the Lebesgue measure on $\T^{n-1}$, \cf \cite[Cor. 3.7]{LugerNedic2017} and \cite[Thm. 4.2]{LugerNedic2021}, see also Lemma \ref{lem:hyperplanes} below. In particular, all points have zero mass, \cf \cite[Rem. 4.3]{LugerNedic2021}. When $n=1$, there are no hyperplanes to consider and point masses are possible.}
    \item{The (topological) support of $\nu$ obeys certain geometric restrictions, \cf \cite[Cors. 4.7 and 4.10]{LugerNedic2021}. When $n=1$, any Borel subset of $\T$ can appear as the support of a RP-measure.}
    \item{If $\nu$ can be written as a product measure, then at least one of the factors must be equal to a constant multiple of the Lebesgue measure, \cf \cite[Prop. 5.11]{Nedic2020}. When $n=1$, there are no product measures to consider.}
    \item{When $n = 2$, the set of probability RP-measures contains an absolutely continuous extremal element, \cf  \cite[pg. 733]{McDonald1990}. Corollary \ref{coro:measure_on_lines} characterizes extremal measures supported on hyperplanes in the case of the poly-upper half-plane. When $n = 1$, the only extremal elements are Dirac measures.}
\end{itemize}

Due to its particular importance for the results of this paper, we present, for the reader's convenience, a proof of the result regarding the restrictions to hyperplanes passing through $\vec{1} \in \T^n$. 

\begin{lemma}\label{lem:hyperplanes}
Let $\nu$ be a RP-measure on $\T^n$. Then, for any hyperplane $M_k:=\{\vec{w} \in \T^n~|~w_k=1\}$ it holds that
$\nu|_{M_k}$ is a multiple of the $(n-1)$-dimensional Hausdorff measure on $M_k$. 
\end{lemma}

\begin{figure}[!ht]
\centering

\begin{tikzpicture}[scale=1.1]
\draw[help lines,->] (-3.2,0) -- (3.2,0) node[above] {$k$};
\draw[help lines,->] (0,-3.2) -- (0,3.2) node[right] {$j_2$};

\fill (0,0)  circle[radius=1.5pt];
\fill (1,0)  circle[radius=1.5pt];
\fill (2,0)  circle[radius=1.5pt];
\fill (3,0)  circle[radius=1.5pt];
\fill (0,1)  circle[radius=1.5pt];
\fill (1,1)  circle[radius=1.5pt];
\fill (2,1)  circle[radius=1.5pt];
\fill (3,1)  circle[radius=1.5pt];
\fill (0,2)  circle[radius=1.5pt];
\fill (1,2)  circle[radius=1.5pt];
\fill (2,2)  circle[radius=1.5pt];
\fill (3,2)  circle[radius=1.5pt];
\fill (0,3)  circle[radius=1.5pt];
\fill (1,3)  circle[radius=1.5pt];
\fill (2,3)  circle[radius=1.5pt];
\fill (3,3)  circle[radius=1.5pt];

\fill (-1,0)  circle[radius=1.5pt];
\fill (-2,0)  circle[radius=1.5pt];
\fill (-3,0)  circle[radius=1.5pt];
\fill (0,-1)  circle[radius=1.5pt];
\fill (-1,-1)  circle[radius=1.5pt];
\fill (-2,-1)  circle[radius=1.5pt];
\fill (-3,-1)  circle[radius=1.5pt];
\fill (0,-2)  circle[radius=1.5pt];
\fill (-1,-2)  circle[radius=1.5pt];
\fill (-2,-2)  circle[radius=1.5pt];
\fill (-3,-2)  circle[radius=1.5pt];
\fill (0,-3)  circle[radius=1.5pt];
\fill (-1,-3)  circle[radius=1.5pt];
\fill (-2,-3)  circle[radius=1.5pt];
\fill (-3,-3)  circle[radius=1.5pt];

\draw [dashed,line width=0.5pt] (-3.2,2) -- (3.2,2);
\node at (-2,2.2) {$\widehat{\nu_2}(\,\cdot\,)$};
\end{tikzpicture}

\caption{Visualization of the possibly non-zero Fourier coefficients of the measures $\nu$ (dots) and $\nu_2$ (dots on the dashed line).}
\label{fig:Fourier_coeff}
\end{figure}
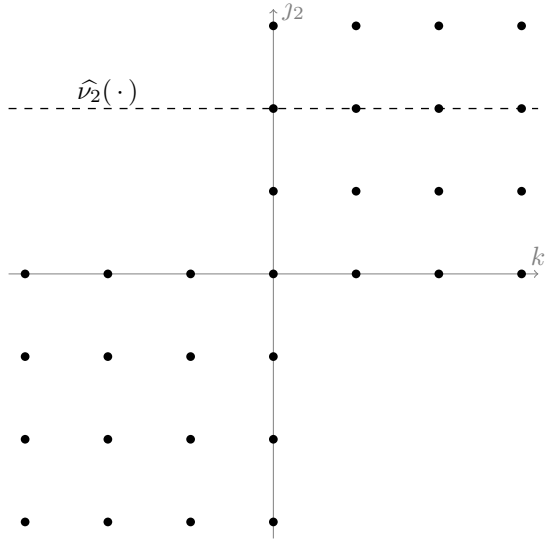

\proof
Without loss of generality, we consider only the hyperplane $M_1$ as all others may be considered analogously. Let us view the measure $\nu|_{M_1}$ as a measure on $\T^{n-1}$ in the variables $w_2,\ldots,w_n$. Take $\vec{j} = (j_2,\ldots,j_n) \in\Z^{n-1}\setminus\{\vec{0}\}$ and define the auxiliary measure $\nu_{\vec{j}}$ on $\T$ for a Borel set $A \subseteq \T$ as
$$\nu_{\vec{j}}(A) := \int_{A \times \T^{n-1}}w_2^{-j_2}\ldots w_n^{-j_n}\diff\nu(\vec{w}).$$
Note now that the multi-index $\vec{j}$ has at least one non-zero entry. Without loss of generality, suppose that $j_2 \neq 0$. Therefore, for the Fourier coefficients of the measure $\nu_{\vec{j}}$, it holds for every $k \in \Z\setminus\{0\}$ with $j_2 \cdot k < 0$ that
$$\widehat{\nu_{\vec{j}}}(k) = \int_{\T}w_1^{-k}\diff\nu_{\vec{j}}(w_1) = \int_{\T^n}w_1^{-k}w_2^{-j_2}\ldots w_n^{-j_n}\diff\nu(\vec{w}) = 0$$
due to the assumption that $\nu$ is a RP-measure, see also Figure \ref{fig:Fourier_coeff}. By the brothers Riesz' theorem \cite[Thm. 17.13]{Rudin1987}, it follows that $\nu_{\vec{j}}$ is absolutely continuous with respect to the Lebesgue measure on $\T$ as either all of its positive or all of its negative Fourier coefficients are zero. In particular, the measure $\nu_{\vec{j}}$ does not have any point masses. Hence,
$$0 = \nu_{\vec{j}}(\{1\}) = \int_{\T^{n-1}}w_2^{-j_2}\ldots w_n^{-j_n}\diff\nu(1,w_2,\ldots,w_n) = \widehat{\nu|_{M_1}}(\vec{j})$$
$\vec{j}\in\Z^{n-1}\setminus\{\vec{0}\}$, implying the desired result by the injectivity of the Fourier transform \cite[Thm. 7.1.5]{HormanderLPDO1}.
\endproof

\begin{remark}
We would like to highlight for the reader the difference in approaches in the proof Lemma \ref{lem:hyperplanes} and \cite[Thm. 4.2]{LugerNedic2021}. The proof of Lemma \ref{lem:hyperplanes} does not rely on any result regarding RP-functions and instead builds upon classical results from measure theory and Fourier analysis. On the other hand, the proof of \cite[Thm. 4.2]{LugerNedic2021} is obtained from the analogous result on the poly-upper half plane, which, in turn, only relies on the integral representation formula \eqref{eq:intRep_Nvar}.
\end{remark}

\subsection{Correspondence between the two classes of functions}
\label{subsec:correspondence}

For the readers sake we recall in this subsection rather carefully, with all details, the explicit correspondence between the Nevanlinna and the RP-measures. This is particularly important for the sake of analysis in Section \ref{sec:fixed_var} below, as there the results for the poly-upper half-plane are deduced by using the ones proven in case of the polydisc.

Denote by $\varphi\colon\C^+ \to \D$ the Cayley transform, \ie
$$\varphi(z) := \frac{z - \I}{z + \I}, \quad z \in \C^+,$$
and let $\varphi^{-1}\colon\D \to \C^+$ be its inverse, \ie
$$\varphi^{-1}(z) := \I\,\frac{1+z}{1-z}, \quad z \in \D.$$
For any \HN function $q$ we may hence define a RP-function $G$ by setting
\begin{equation}
    \label{eq:q_to_G}
    G(\vec{z}):= -\I\,q\big(\varphi^{-1}(z_1),\ldots,\varphi^{-1}(z_n)\big), \quad \vec{z} \in \D^n.
\end{equation}
Conversely, starting with any RP-function $G$, we may define a \HN function $q$ by setting
\begin{equation}
    \label{eq:G_to_q}
    q(\vec{z}) := \I\,G\big(\varphi(z_1),\ldots,\varphi(z_n)\big),\quad \vec{z} \in \C^{+n}.
\end{equation}
For any two function $q$ and $G$ connected via relations \eqref{eq:q_to_G} and \eqref{eq:G_to_q}, one establishes the following relations between the representing parameters $(a,\vec{b},\mu)$ of $q$ and $(A,\nu)$ of $G$. For the constant factors, it holds that
$$A = \Im[G(\vec{0})] = \Re[q(\I\,\vec{1})] =a.$$
Define now $\beta_k,\,k = 1,\ldots,n$, to be the measure on $\T^n$ where we take the normalized Lebesgue measure $\lambda_\T$ in each variable except the $k$-th, where we take the Dirac measure $\delta_1$ instead, \ie
$$\beta_k := \lambda_\T \otimes \ldots \otimes \lambda_\T \otimes \delta_1 \otimes \lambda_\T \otimes \ldots \otimes \lambda_\T.$$
We note that $\beta_k$ is, in fact, a RP-measure and it holds that the \HN function $q(\vec{z}) = b_kz_k,\,b_k \geq 0,$ corresponds to the RP-function $G(\vec{z}) = b_k\,\frac{1+z_k}{1-z_k}$ whose representing measure equals $b_k\,\beta_k$. Conversely, since we know that the restriction of any RP-measure $\nu$ to a hyperplane passing through $\vec{1} \in \T^n$ is a constant multiple of the Lebesgue measure, we may write $\nu|_{M_K} = d_k\beta_k$ for some $d_k \geq 0$, allowing us to also validate the above correspondence in the reverse direction.

Finally, let $\Phi \colon \R^n \to (0,2\pi)^n$ be defined as
$$\Phi(\vec{t}) := \big(-\I\,\Log(\varphi(t_1)),\ldots,-\I\,\Log(\varphi(t_n))\big),$$
where the logarithm $\Log$ is take to have a branch cut along the positive real axis. Its inverse $\Phi^{-1} \colon (0,2\pi)^n \to \R^n$ is then defined as
$$\Phi^{-1}(\vec{s}) := \big(\varphi^{-1}(\E^{\I\,s_1}),\ldots,\varphi^{-1}(\E^{\I\,s_n})\big),$$
These maps allows us to hop back and forth between measures on $\R^n$ and measures on $(0,2\pi)^n$ and amount to using the change of variables $\E^{\I\,s} = \varphi(t)$ between $s \in (0,2\pi)$ and $t \in \R$ in each coordinate separately. Furthermore, we note that any measure on $(0,2\pi)^n$ that is viewed as measure on $[0,2\pi)^n$ is assumed to be extended by zero on $[0,2\pi)^n\setminus(0,2\pi)^n$. Out of any Nevanlinna measure $\mu$ on $\R^n$ we may thus define a measure $\nu$ on $\T^n$ by setting
\begin{equation}
    \label{eq:mu_to_nu}
    \nu := \Sigma^{-1}(\mu \circ \Phi^{-1}),
\end{equation}
which will be, in fact, a RP-measure. Conversely, starting with a RP-measure $\nu$, we may define a measure $\mu$ by setting
\begin{equation}
    \label{eq:nu_to_mu}
    \mu := \Sigma(\nu|_{(\T\setminus\{1\})^n}) \circ \Phi,
\end{equation}
which will be a Nevanlinna measure. Note also that transformations \eqref{eq:mu_to_nu} and \eqref{eq:nu_to_mu} preserve the singular continuous and absolutely continuous parts of the measures, \ie
$$(\Sigma^{-1}(\mu \circ \Phi^{-1}))_\mathrm{a.c.} = \Sigma^{-1}(\mu_{\mathrm{a.c.}} \circ \Phi^{-1}),$$
etc. When $n=1$, point masses are also preserved.

To summarize, given a \HN function $q$ with representing parameters $(a,\vec{b},\mu)$, the function $G$ given by formula \eqref{eq:q_to_G} will be a RP-function with representing parameters
\begin{equation}
    \label{eq:q_to_G_parameteres}
    \big(a,\textstyle{\sum_{k=1}^n}b_k\beta_k + \Sigma^{-1}(\mu \circ \Phi^{-1})\big)
\end{equation}
Conversely, given a RP-function $G$ with representing parameters $(A,\nu)$, the function $q$ given by formula \eqref{eq:G_to_q} will be a \HN function with representing parameters
\begin{equation}
    \label{eq:G_to_q_parameters}
    \big(A,(d_1,\ldots,d_n),\Sigma(\nu|_{(\T\setminus\{1\})^n}) \circ \Phi\big),
\end{equation}
where $d_k \geq 0$ are the constants for which $\nu|_{M_K} = d_k\beta_k$.

Finally, we wish to emphasize again that the validity of transformations \eqref{eq:mu_to_nu} and \eqref{eq:nu_to_mu} between Nevanlinna and RP-measures, and consequentially relations \eqref{eq:q_to_G_parameteres} and \eqref{eq:G_to_q_parameters}, is dependant on fact that the map $\Psi$ and relations \eqref{eq:nu_to_sigma}, and consequentially relations \eqref{eq:q_to_G} and \eqref{eq:G_to_q}, are defined exactly as specified. One can, of course, always pick a different parametrization of $\T^n$ or a different biholomorphism between $\C^+$ and $\D$. However, in such a case, all formulas must be re-adjusted accordingly. 

\section{Characterization of Nevanlinna measures via the Fourier transform}
\label{sec:Fourier}

Recall from Section \ref{subsec:RP_functions} that the RP-measures on the polydisc $\T^n$ are characterized simply as the measures whose non-zero Fourier coefficients are supported on the union of the positive and negative 'octants', \ie if $\nu(\vec{j}) \neq 0$ for $\vec{j} \in \Z^n$, then either $j_\ell \geq 0$ for all $1\leq \ell \leq n$ or, alternatively, $j_\ell \leq 0$ for all $1 \leq \ell \leq n$. A natural question arises whether the same is true for Nevanlinna measures on $\R^n$, when we replace the Fourier series by the Fourier transform. The answer turns out to be positive, see Theorem \ref{thm:Fourier_char} below. It is somewhat surprising that this fact seems not to have been noticed in the literature before. One reason for this could be that in the generic case the Fourier transform of a Nevanlinna measure will not be a function (or a measure), but has to be understood as a Schwartz distribution. The formal idea of the proof will be rather straightforward, but one needs to work quite a bit in order to ensure that the formal computations make sense.

\subsection{Definitions and prerequisites}

Let us start by assuming that $\mu$ is a positive Borel measure on $\R^n$ satisfying the necessary growth condition \eqref{eq:measure_growth}, which we rephrase for further purposes as a finiteneness of a norm:
\begin{equation}\label{eq:001}
\|\mu\|_{\NL}:=\int_{\R^n}\prod_{\ell=1}^{n}\frac{1}{1+t_\ell^2}\diff|\mu|(\vec{t}) < \infty.
\end{equation}
We denote by $\NL(\R^n)$ the Banach space of all signed measures $\mu$ on $\R^n$ for which he norm $\|\mu\|_{\NL}$  is finite, and by $\NLP(\R^n)$ the cone of non-negative measures in $\NL(\R^n)$. Then, the Poisson extension of $\mu$ is well-defined for $(\vec{x},\vec{y})\in \R^n \times \R_+^n$ as 
$$
u(\vec{x},\vec{y}) = \frac{1}{\pi^n}\int_{\R^n}\pois_n(\vec{x}+\I\,\vec{y},\vec{t})\diff\mu(\vec{t}) = \frac{1}{\pi^n}\int_{\R^n}\prod_{\ell=1}^n\frac{y_\ell}{(x_\ell-t_\ell)^2+y_\ell^2}\diff\mu(\vec{t}),
$$
where $\pois_n$ denotes the Poisson kernel of $\C^{+n}$ as defined in formula \eqref{eq:Poisson_half_plane}. Note that the extension $u$ is separately harmonic\footnote{It is due to this fact that we have chosen the notation $\NL$.} in $\C^{+n}$ and will be pluriharmonic in $\C^{+n}$ if and only if the measure $\mu$ also satisfies the Nevanlinna condition \eqref{eq:measure_Nevan} \cite[Prop. 5.2]{LugerNedic2019}. Because of the mild growth ensured by condition \eqref{eq:001}, the measure $\mu$ defines a Schwartz distribution on $\R^n$. The same is true for the restrictions $\vec{x} \mapsto u(\vec{x},\vec{y})$ for a fixed point $\vec{y} \in \R_+^n$, or for the derivatives $\vec{x} \mapsto D^{\vec{\alpha}}_{\vec{y}}u(\vec{x},\vec{y})$.

We need to express the quantities we want to study in terms of the Fourier transform $\widehat u(\vec{\xi},\vec{y}),$ where the Fourier transform will always be taken with respect to the $\vec{x}$-variables. Recall that the Fourier transform $\four$ of a function $f \in \Leb^1(\R^n)$ is defined as
$$(\four f)(\vec{\xi}) := \widehat{f}(\vec{\xi}) = \int_{\R^n}f(\vec{x})\E^{-\I\,\vec{x}\cdot\vec{\xi}}\diff\vec{x},$$
where $\vec{x}\cdot\vec{\xi} := x_1\xi_1 + \ldots + x_n\xi_n$. Then the inverse Fourier transform takes the form $f(x) = (2\pi)^{-n}\int_{\R^n}(\four f)(\vec{\xi})\E^{\I\,\vec{x}\cdot\vec{\xi}}\diff\vec{\xi}.$

In order to get us moving, let us first assume that $\mu(\R^n)<\infty$. Then, by recalling the one-dimensional Fourier transform $\big(\four(1+x^2)^{-1}\big)(\xi)= \pi \E^{-|\xi|}$, we obtain that the Fourier transform of the Poisson extension of $\mu$ equals
\begin{equation}\label{eq:003}
\widehat u(\vec{\xi},\vec{y}) = \prod_{\ell=1}^n \E^{-y_\ell|\xi_\ell|}\,\cdot\,\widehat\mu(\vec{\xi}),
\end{equation}
where the product is well-defined since $\widehat\mu$ is continuous in this case. In a similar fashion, by differentiating above with respect to $y_k$ and $x_k$, we see that 
\begin{equation}\label{eq:005}
\four\big(\partial_{z_k}\partial_{\overline{z_j}}u)(\vec{\xi},\vec{y}) = m_{\vec{y},k,j}(\vec{\xi})\widehat\mu(\vec{\xi})\
\end{equation}
where $\partial_{z_k} = (\frac{\diff}{\diff x_k}-\I\,\frac{\diff}{\diff y_k})/2$, \, $\partial_{\overline{z_k}} = (\frac{\diff}{\diff x_k} + \I\,\frac{\diff}{\diff y_k})/2$ and
\begin{equation}
    \label{eq:Fourier_multiplier}
    m_{\vec{y},k,j}(\vec{\xi}) := -\,\xi_k\chi_{\{ \xi_k<0\}}(\vec{\xi})\cdot\xi_j\chi_{\{ \xi_j>0\}}(\vec{\xi}) \cdot \prod_{\ell=1}^n \E^{-y_\ell|\xi_\ell|}.
\end{equation}

In order to deal with the case where $\mu$ is not a finite measure, we need to consider a suitable class of Fourier multipliers. Let $m \colon\R^n \to \C$  be a bounded continuous function. Then the action $T_m\mu ={\mathcal F}^{-1}(m\widehat f)$ is well-defined  and an element in $S'(\R^n)$ in case  $\mu$ is a finite measure, since then $m\widehat f$ is a bounded and continuous function. The next lemma gives a sufficient condition ensuring that $T_m$ extends to all of $\NL.$  We denote by $\|\mu\|_{TV}$ the total variation norm of a signed measure.

\begin{lemma}\label{lem:Fourier_multiplier}
Assume that $m\colon\R^n\to\C$ is a bounded continuous function such that 
$$
\|D^{\vec{\alpha}} m\|_{TV} <\infty,
$$
for all multi-indices $\vec{\alpha} \in \N_0^n$ such that $\alpha_k\leq 2$ for all $k=1,\ldots, n$. Then, $m$ defines a  multiplier on $\NL$ and the action $T_m\mu$ is  well-defined by setting
\begin{equation}\label{eq:019}
T_m\mu := \lim_{\ell\to\infty}T_m\mu_\ell,
\end{equation}
where $\mu_\ell$:s are finite measures with $\|\mu_\ell-\mu\|_{\NL}\to 0$ as $\ell\to\infty.$ Moreover, the multiplier $T_m$ is continuous on $\NL$, i.e.  if $\|\mu_\ell-\mu\|\to 0$ as $\ell\to\infty$, where $\mu, \mu_1,\mu_2,\ldots$ are elements if $\NL$, then
$$
T_m \mu_\ell\to T_m\mu \qquad\textrm{in}\quad \schw'(\R^n).
$$
Actually, $T_m\mu$ is continuous and the convergence above is uniform on compact sets.
\end{lemma}

\proof
It is enough to prove that for any $R>0$  there is a constant $C<\infty$ such that for all finite signed measures $\mu$  we have
\begin{equation}\label{eq:009}
\sup_{|\vec{x}|\leq R} |(T_m\mu)(\vec{x})|\leq C\|\mu\|_{\NL}(R+1)^{2n}.
\end{equation}
Namely, this implies that the definition \eqref{eq:019} is well-posed since using the approximation $\mu_k:= \chi_{B(0,k)}\mu$ we infer by \eqref{eq:009} that $\four(T_m\mu)$  is continuous with at most polynomial growth, and hence defines an element $S'(\R^n)$. Moreover, if $\mu_\ell,\mu$, $\ell\geq 1$, are arbitrary measures in $\NL$ and such that  $\lim_{\ell\to\infty}\|\mu_\ell-\mu\|_{\NL}= 0$, we deduce that  $(1+|\vec{x}|^2)^{-n}T_m\mu_\ell \to (1+|\vec{x}|^2)^{-n}T_m\mu$ uniformly on $\R^n$ which yields the stated convergence in the sense of distributions, and also the unifrom convergence on compact sets.

In order to prove \eqref{eq:009}, define the finite measure $\nu$ by setting 
$$
\diff\nu(\vec{t})= \prod_{\ell=1}^{n}\frac{1}{1+t_\ell^2}\diff\mu(\vec{t})
$$
and assume that $|x|\leq R$. We obtain
\begin{eqnarray*}
T_m\mu(\vec{x})& = & \frac{1}{(2\pi)^n}\int_{\R^n} \E^{\I\,\vec{\xi}\cdot\vec{x}}m(\vec{\xi})\widehat\mu(\vec{\xi})\diff\vec{\xi} \\
~ & = & \frac{1}{(2\pi)^n}\int_{\R^n}\E^{\I\,\vec{\xi}\cdot \vec{x}}m(\vec{\xi})\big(\prod_{\ell=1}^{n}(1-(\diff/\diff\xi_\ell)^2)\big)\widehat\nu(\vec{\xi}) \diff\vec{\xi} \\
&=& \frac{1}{(2\pi)^n}\int_{\R^n} \Big[\big(\prod_{\ell=1}^{n}(1-(\diff/\diff\xi_\ell)^2)\big)\E^{\I\,\vec{\xi} \cdot \vec{x}}m(\vec{\xi})\Big]\widehat\nu(\vec{\xi})\diff\vec{\xi} \\
&=& \int_{\R^n} \widehat\nu(\vec{\xi}) \diff\eta_{\vec{x}}(\vec{\xi}),
\end{eqnarray*}
where $\eta_{\vec{x}}$ is the signed measure
$$
\eta_{\vec{x}} := (2\pi)^{-n} \cdot \big(\prod_{\ell=1}^{n}(1-(\diff/\diff\xi_\ell)^2)\big)\E^{\I\,\vec{\xi} \cdot \vec{x}}m(\vec{\xi}).
$$
Our assumption on the multiplier $m$ verifies that 
$$
\|\eta_{\vec{x}}\|_{TV}\leq C(1+|\vec{x}|^2)^n,
$$
and the claim follows by combining this with 
$$
\|\widehat\nu\|_{L^\infty(\R^n)}\leq \nu(\R^n)=\|\mu\|_{\NL}<\infty.
$$
This finishes the proof.
\endproof

\subsection{Characterization theorem and first corollaries} 

The following theorem  gives the announced characterization of Nevanlinna measures in terms of the support of their Fourier transform.

\begin{thm}\label{thm:Fourier_char}
Assume that $\mu$ is a positive Borel measure on $\R^n$ with $n \geq 2$ satisfying the growth condition \eqref{eq:measure_growth}. Then, $\mu$ is a Nevanlinna measure if and only its distributional Fourier transform satisfies
$$
\supp(\widehat\mu) \; \subset\; [0,\infty)^n\cup  (\infty,0]^n.
$$
\end{thm}
\proof
Assume first that $\mu$ is a Nevanlinna  measure. Choose a non-negative, symmetric, radial, and smooth test-function $g\colon\R^n \to \R$ with support in $B(\vec{0},1)\subset\R^n$ and with $\int_{\R^n}g(\vec{x})\diff\vec{x}=1$. For $\varepsilon >0$, denote $G_\varepsilon(\vec{x}) := \widehat g(\varepsilon \vec{x})$ and $g_\varepsilon (\vec{\xi}):=\varepsilon^{-n}g(\vec{\xi}/\varepsilon)$. Note that $G_\varepsilon$ decays at any polynomial rate as a Schwartz function, and we have that 
$$
\mu_\varepsilon:=G_\varepsilon\mu
$$
is a finite measure with smooth Fourier transform $\widehat\mu_\varepsilon=g_\varepsilon*\widehat\mu.$ Since $\|\mu_\varepsilon-\mu\|_{\NL}\to 0$ as $\varepsilon\to 0$, we have that $\pois_n\mu_\varepsilon \to \pois_n\mu$ locally uniformly in $\C^{+n}$, and there is also convergence in $\schw'(\R^n)$. Moreover, for any fixed $\vec{y}\in\R_+^n$, we have
\begin{equation}\label{eq:011}
\pois_n\mu_\varepsilon(\,\cdot\,,\vec{y}) \to \pois_n\mu(\,\cdot\,,\vec{y}) \quad \textrm{in}\;\; \schw'(\R^n)\quad\textrm{as}\quad \varepsilon\to 0,
\end{equation}
and then the convergence in $\schw'(\R^n)$ also follows for the Fourier transforms of the extension. The Fourier multiplier $m_{\vec{y},k,j}$ defined in formula \eqref{eq:Fourier_multiplier} is admissible by Lemma \ref{lem:Fourier_multiplier}. Hence we may compute, for any $k,j\in\{1,\ldots, n\}$ with $k \neq j$, that
\begin{eqnarray}\label{eq:013}
0 = \four\big(\partial_{z_k}\partial_{\overline{z_j}}u(\,\cdot\,,\vec{y})\big)(\vec{\xi}) = \lim_{\varepsilon \to 0} m_{\vec{y},k,j}(\vec{\xi})\widehat\mu_\varepsilon(\vec{\xi}),
\end{eqnarray}
with convergence in $\schw'(\R^n)$. This especially implies that $\widehat\mu$ vanishes in any open set where $m_{\vec{y},k,j}(\vec{\xi})$ is non-zero and smooth, \ie on the 'quadrant'
$$\{\vec{\xi} \in \R^n~|~ \xi_k < 0 \;\text{and}\; \xi_j>0\}.$$
Since this holds for any $k \neq j$, the claim on the support of $\widehat \mu$ follows.

Towards the other direction, let us assume that $\widehat\mu$ has the stated support and fix $k,j\in\{1,\ldots, n\}$, $k \neq j$, together with $\vec{y} \in \R_+^n$. We need to show that $\partial_{z_k}\partial_{\overline{z_j}}u=0$. We cannot deduce it directly from \eqref{eq:013} since $m_{\vec{y},k,j}(\xi)$ is not vanishing in a full neighbourhood of support of $\widehat\mu$, and also the non-smoothness of $m_{\vec{y},k,j}(\vec{\xi})$ could be problematic. Thus let us fix $\delta >0$ and consider the translated multiplier
$$
m_{\vec{y},k,j,\delta}(\vec{\xi}) := m_{\vec{y},k,j}(\vec{\xi}+\delta(\vec{e}_k-\vec{e}_j)),
$$
where $\vec{e}_\ell$ denotes the $\ell$-th standard basis vector of $\R^n$. Note that this multiplier vanishes in the set $\supp(\mu)+\varepsilon B(\vec{0},1)$ as soon as $\varepsilon < \delta$. Thus, we see that
\begin{equation}\label{eq:015}
T_{m_{\vec{y},k,j,\delta}}\mu = \lim_{\varepsilon \to 0}T_{m_{\vec{y},k,j,\delta}} \mu_\varepsilon =0.
\end{equation}
The desired conclusion $\partial_{z_k}\partial_{\overline{z_j}}u(\,\cdot\,,\vec{y})=0$ follows from \eqref{eq:015} by noting that
\begin{eqnarray}\label{eq:021}
\partial_{z_k}\partial_{\overline{z_j}}u(\,\cdot\,,\vec{y}) = T_{m_{\vec{y},k,j}}\mu = \lim_{\delta\to 0^+} T_{m_{\vec{y},k,j,\delta}}\mu
\end{eqnarray}
since we may compute the action of a translate of a multiplier as follows:
$$
T_{m(\,\cdot\,+ \delta \vec{a})}\mu (\vec{x})= \E^{-\delta \vec{a} \cdot \vec{x}} T_m \big( \E^{\delta \vec{a} \cdot \vec{x}}\mu\big)(x),
$$
where $\vec{a} \in \R^n$, and then \eqref{eq:021} follows by Lemma \ref{lem:Fourier_multiplier} together with the observation  
$$
\|\E^{\delta \vec{a} \cdot \vec{x}}\mu-\mu\|_{\NL}\to 0\quad \textrm{as}\quad \delta\to 0.
$$
This finishes the proof.
\endproof

\begin{remark}\label{rem:one_dimensional_Fourier}
One should note that in dimension $n=1$, Theorem \ref{thm:Fourier_char} puts no restriction on the measure $\mu$ besides the necessary condition $\|\mu\|_{\NL}<\infty$, which, of course, coincides with the standard results from the one-dimensional theory.
\end{remark}

\begin{remark}\label{rem:hormander} One direction of the above theorem could also be deduced from extended versions of Paley-Wiener theorem, see \eg \cite[Sec. 7.3]{HormanderLPDO1} or \cite[Sec. 2.2]{Simon1974}.
\end{remark}

\begin{example}\label{ex:Fourier_examples_dim2}
Let us consider some examples in dimension $n=2.$ As will be checked in connection with Corollary \ref{coro:measure_on_lines} below, if $L$ is a line through origin with a negative slope, and $\mu$ is the restriction of the 1-dimensional Hausdorff measure on $L$, then $\widehat\mu$ is supported on a line $L'$ perpendicular to $L$, and hence $\mu$ is a Nevanlinna measure according to Theorem \ref{thm:Fourier_char}. On the other hand this can be also be seen by noting that up to a constant, $\mu$ is the boundary distribution corresponding to the function
$$(z_1,z_2) \mapsto \Im[(-a(z_1+b z_2)^{-1}]$$
for suitable constants $a,b > 0$.

For our second example, let us first recall the Fourier transform of a power function on the positive real axis. Assume that  $\alpha \in (-1,1)$ and $\varepsilon, u>0$. We get
$$
\int_0^\infty x^\alpha \E^{-(\varepsilon+u)x}\diff x = \Gamma (\alpha+1)(\varepsilon +u)^{-1-\alpha}
$$
where $\Gamma$ denotes Euler's Gamma functions. By analytic continuation in $u$ in the region $\Re[u] \geq 0$ we may set $u=\I\,\xi$ and deduce that
$$\four(x_+^\alpha \E^{-\varepsilon\,x})(\xi) = \Gamma (\alpha+1)(\varepsilon + \I\,\xi)^{-1-\alpha},$$
where $x_+ := \max\{0,x\}$ and the branch of $z^{-\alpha-1}$ (considered in $ \Re[z] \geq \varepsilon$) is real on the positive real axis. Letting $\varepsilon\to 0^+$ we obtain that
\begin{equation}\label{eq:039}
\four(x_\pm^\alpha) (\xi)= \E^{\mp\,\I\,\pi\,(\alpha+1)\, \sgn(\xi)/2}\Gamma (\alpha+1)|\xi|^{-1-\alpha}, \quad |\xi|>0,
\end{equation}
where $x_- := \max\{0,-x\}$. At the origin, the Fourier transform may be defined as the limit distribution  as $\varepsilon\to 0^+$.

Consider now the measure $\mu_\alpha$ on $\R^2$, where $\mu_\alpha$ has the density
\begin{multline}\label{eq:043}
d\mu_\alpha(x) \\ = \big(a_{11}(x_1)_+^\alpha(x_2)_+^\alpha +a_{12} (x_1)_+^\alpha(x_2)_-^\alpha+a_{21} (x_1)_-^\alpha(x_2)_+^\alpha+a_{22}(x_1)_-^\alpha(x_2)_-^\alpha \big)\diff x_1 \diff x_2,
\end{multline}
where $a_{11},a_{12},a_{21},a_{22} \geq 0$. We will use Theorem \ref{thm:Fourier_char} to check what relations the coefficients have to satisfy so that $d\mu_\alpha$ is a Nevanlinna measure in $\R^2$. To that end, we use formula  \eqref{eq:039} in order to compute that 
\begin{eqnarray*}
\widehat\mu_\alpha(\xi)&=&
A(\Gamma (\alpha+1))^2|\xi_1\xi_2|^{-1-\alpha} \quad\textrm{for}\quad \xi_1<0, \; \xi_2>0\qquad\textrm{and} \\
\widehat\mu_\alpha(\xi)&=&
B(\Gamma (\alpha+1))^2|\xi_1\xi_2|^{-1-\alpha} \quad\textrm{for}\quad \xi_1>0, \; \xi_2<0,
\end{eqnarray*}
where 
\begin{eqnarray*}\label{eq:045}
A &:=& a_{11}+a_{12}\,\E^{\I\,\pi\,(\alpha+1)}+a_{21}\,\E^{-\I\,\pi\,(\alpha+1)}+a_{22}\,\qquad \textrm{and}\\
B &:=& a_{11}+a_{12}\,\E^{-\I\,\pi\,(\alpha+1)}+a_{21}\,\E^{\I\,\pi\,(\alpha+1)}+a_{22}.
\end{eqnarray*}
Thus, the necessary and sufficient condition for $\mu$ being Nevanlinna is given by the requirement $A=B=0$, which due to $A=\overline{B}$ reduces to $A=0$. The imaginary part of $A$ is zero if $a_{21}=a_{12}$, or if $\alpha=0$, and then the real part vanishes only if 
\begin{eqnarray}\label{eq:047}
a_{11}+a_{22}=\cos(\pi\alpha)(a_{12}+a_{21}).
\end{eqnarray}
This clearly implies that a necessary condition is $|\alpha|\leq 1/2$. In case $\alpha=\pm 1/2$ we see that the only solution is $a_{11}=a_{22}=0$ and $a_{12}=a_{21}$.
In case $|\alpha|<1/2$ the necessary and sufficient conditions (beside non-negativity) are given by $a_{12}=a_{21}$ and relation \eqref{eq:047}. 

Finally, note that an explicit example of a \HN function with such a representing measure is the function $(z_1,z_2) \mapsto \sqrt{z_1}\sqrt{z_2}$, where the branch cut of the square root is taken along the negative real axis. The representing measure $\mu$ of this function has the density
$$\diff\mu = \big((x_1)_-^{\frac{1}{2}}(x_2)_+^{\frac{1}{2}} + (x_1)_+^{\frac{1}{2}}(x_2)_-^{\frac{1}{2}}\big)\diff x_1 \diff x_2,$$
\ie $a_{11} = a_{22} = 0$ and $a_{12} = a_{21} = 1$ with $\alpha = \frac{1}{2}$.\hfill$\lozenge$
\end{example}

Next, we present some applications of Theorem \ref{thm:Fourier_char}. 

\begin{coro}[\cf {\cite[Prop. 4.3]{LugerNedic2017}}]
\label{coro:finite_measure}
If $\mu$ is a non-zero Nevanlinna measure in $\R^n$ with $n \geq 2$, then $\mu(\R^n)=\infty$.
\end{coro}

\proof
If $\mu$ is a finite positive Borel measure, its Fourier transform is continuous with $\widehat \mu(\vec{0})=\mu(\R^n) > 0$, whence $\widehat \mu$ cannot vanish in any octant, and hence $\widehat \mu$ does not satisfy the support condition of Theorem \ref{thm:Fourier_char}.
\endproof

\begin{coro}[\cf {\cite[Thm. 4.1]{Nedic2020}}]
\label{coro:product_measure}
Let $\mu_1$ be a positive measure in $\R^{n_1}$ and $\mu_2$ be a positive measure in $\R^{n_2}$. Then, the product measure $\mu_1\times \mu_2$ is a non-zero Nevanlinna measure in $\R^{n_1+n_2}$ if and only if one of the factors is a multiple of the Lebesgue measure and the other one is a Nevanlinna measure in the respective dimension. 
\end{coro}

\proof
Just use Theorem \ref{thm:Fourier_char} and note that the factors of the product $\widehat{\mu\times\mu_2}=\widehat {\mu_1}\times \widehat {\mu_2}$ have symmetric supports with respect of the origin. Hence the product cannot have support in $[0,\infty)^n\cup  (\infty,0]^n$ unless one of the factors is supported at the origin, and hence is a finite sum of the derivates of $\delta_0$. Moreover, this implies that the corresponding factor has a polynomial density with respect to the Lebesgue measure, and the density is easily seen to be constant in view of \eqref{eq:measure_growth}. The converse implication also follows immediately from  \ref{thm:Fourier_char}.
\endproof

\subsection{Measures on hyperplanes}

The following theorem characterizes Nevanlinna measures supported on an arbitrary hyperplane in $\R^n$. Note that the hyperlanes contained in a translate of one of the $(n-1)$-dimensional coordinate hyperplanes are already taken care by Corollary \ref{coro:product_measure}.

\begin{thm}\label{thm:hyperplanes}
 Let $M\subset \R^n$ be a $m$-dimensional hyperplane that is not parallel to one of the coordinate hyperplanes. Then $M$ can support a non-trivial Nevanlinna measure if and only if $m=n-1$ and there is $\vec{a}=(a_1,\ldots, a_n) \in\R^n$ with $a_k\geq 0$
for all $k$, and at least two of the $a_k$-s non-zero, and $c\in\R$ such that
$$
M=\{ \vec{x}\in \R^n~|~\vec{a}\cdot \vec{x}=c\}.
$$
In this case, the Nevanlinna measures supported by $M$ are then exactly the measures
$$
\mu= p(x_1,\ldots ,x_n)\mathcal{H}^{m}_{M},
$$
where the density $p$ is a polynomial of at most degree 2 such that $p_{|M}\geq 0$, and $p$ depends only on the variables $x_j$ for which $a_j>0$.
\end{thm}

\proof
By translation it is enough to consider a $m$-dimensional hyperplane $M\subset\R^n$ through the origin, where $1\leq m\leq n-1$, not contained in any of the hyperplanes $\{x\in\R^n \; |\; x_j=0\}$, $1\leq j\leq n$. Assume that $\mu$ is a Nevanlinna measure carried by $M$. Thus, we may write  in the natural coordinates $\mu=\nu\times \delta_0$,  where  $\nu$ is a positive measure on $M$, and $\delta_0$ is the delta measure at the origin on $M^\perp$.   Then, again in the natural coordinates, we have 
$$
\widehat \mu =\widehat\nu\times (\mathcal{H}^{n-m}_{|M^\perp}),
$$
where  $\mathcal{H}^{n-m}_{|M^\perp}$ is the restriction of the $(n-m)$-dimensional Hausdorff measure to $M^\perp$. Theorem \ref{thm:Fourier_char} implies that  $\mu$ is a Nevanlinna measure if and only if it satisfies growth \eqref{eq:measure_growth} together with the inclusion
\begin{equation}\label{eq:perus}
\supp(\widehat \mu) = \supp(\widehat\nu)\times M^\perp\subset (-\infty,0]^n\cup [0,\infty)^n =:C_-\cup C_+,
\end{equation}
where  naturally ${\rm supp\,}(\widehat \nu)$ is considered a subset of $M$. Since $\widehat \nu$ is the Fourier transform of a positive measure, its support contains the origin of $M$. Hence \eqref{eq:perus} especially implies that
$$
M^\perp\subset (-\infty,0]^n\cup [0,\infty)^n.
$$
However, this forces $M^\perp$ to be one-dimensional. Namely, if this would not be true, we could choose two non-zero vectors $\vec{b},\vec{b'}\in M^\perp$ such that  $\vec{b}\cdot\vec{b'}=\sum_{j=1}^nb_jb'_j=0$. Then either both $\vec{b}$ or $\vec{b'}$ have coordinates with varying signs, or otherwise we have $\sum_{j=1}^n|b_jb'_j|=0$. In any case, one of the vectors $\vec{b}, \vec{b}\pm \vec{b'}$ has coordinates with varying signs, and they cannot all be contained in $C_-\cup C_+$, which yields  a contradiction. Hence dim($M^\perp)=1$ and  we may pick a unit vector $\vec{a} \in C_+$ such that $M^\perp=\{t\,\vec{a}~|~t\in\R\}$ so that
$$
M=\{ x\in\R^n~|~\vec{a} \cdot \vec{x}\}.
$$
Also, at least two coordinates of $\vec{a}$ are non-zero, since otherwise $M$ would be a coordinate hyperplane.

Next, we claim that $\supp(\widehat\nu)=\{0\}$.  Assume that this is not true, and there is $\vec{b}\in M$ such that $\vec{b}\in \supp(\nu)\setminus\{0\}$. Then \eqref{eq:perus} implies that $\vec{b}+M =  \{\vec{b}+t\,\vec{a}~|~t\in\R\}\subset C_-\cup C_+$. Since $\vec{b}$ and $\vec{a}$ are orthogonal, we deduce as before that one of the vectors $\vec{b}, \vec{b}\pm \vec{a}$ is not in $C_-\cup C_+$, which contradicts \eqref{eq:perus}. Thus $\widehat{\nu}$ has to be supported at the origin of $M$, which shows that it is a finite sum of derivatives of $\delta$. Then,  in any orthogonal coordinates  on $M$, the measure $\nu$ has a polynomial density, which naturally is also a polynomial density on $M$ in terms of the full set of coordinates $x_1, \ldots , x_n$. This density naturally needs to satisfy the  growth condition \eqref{eq:measure_growth}.

Conversely, if  $M$ is of the stated form, and $\nu$ is any polynomial density on $M$ satisfying \eqref{eq:measure_growth}, then it is also a polynomial density on $M$ with respect to arbitarily chosen orthornnormal coordinates on $M$, and it follows that then
$\supp(\widehat \mu)= M^\perp=\{t\,\vec{a}~|~t\in\R\}\subset C_-\cup C_+$  , whence $\mu$ is a Nevanlinna measure by Theorem \ref{thm:Fourier_char}.

Finally, we need to check which polynomials densities satisfy condition \eqref{eq:measure_growth}. To that end, by symmetry, we may assume that $a_1,\ldots, a_\ell>0$ and $a_{\ell+1},\ldots a_n=0$ for some $2\leq \ell\leq n$. In this case we may parametrize $M$ by the coordinates $x_1,x_2,\ldots x_\ell-1, x_{\ell+1},\ldots, x_n$ via the bijective map
$$
\R^{n-1}\ni (x_1,x_2,\ldots x_{\ell-1}, x_{\ell+1},\ldots, x_n)\mapsto (x_1,x_2,\ldots x_{\ell-1}, \sum_{j=1}^{\ell-1}b_jx_j, x_{\ell+1},\ldots, x_n)\in M,
$$
where $b_j=-a_j/a_\ell\not=0$. This parametrization preserves the $\mathcal{H}^{n-1}$-measure up to a multiplicative non-zero constant. Any polynomial density can be written as a polynomial in these coordinates of $M$, and it remains to check which positive polynomials $p$ satisfy the condition
$$
\int_{\R^{n-1}}\frac{p(x_1,\ldots, x_{\ell-1},x_{\ell+1},\ldots, x_n)\diff x_1\ldots \diff x_{\ell-1} \diff x_{\ell+1}\ldots \diff x_n}
{\big(1+x_1^2\big)\ldots \big(1+x_{\ell-1}^2\big)\big(1+(\sum_{j=1}^{\ell-1}b_jx_j)^2\big)\big(1+x_{\ell+1}^2\big)\ldots \big(1+x_n^2\big)} \; <\infty.
$$
If the polynomial $p$ has nontrivial dependence on any one of the variables $x_{\ell+1}, \ldots, x_n$, say $x_k$ with $\ell+1\leq k\leq n$, then for almost every value of the other variables the growth of $|p|$ in $x_{k}$ is at least linear as $|x_k|\to\infty$, and hence an application of Fubini's theorem shows that the above integral is not finite. Thus, the polynomial $p$ may only depend on the variables $x_1,\ldots, x_{\ell-1}$, and we are reduced to finding non-negative polynomials $p=p(x_1,\ldots, x_{\ell-1})$ such that 
 $$
\int_{\R^{\ell-1}}\frac{p(x_1,\ldots, x_{\ell-1})dx_1\ldots dx_{\ell-1} }
{\big(1+x_1^2\big)\ldots \big(1+x_{\ell-1}^2\big)\big(1+(\sum_{j=1}^{\ell-1}b_jx_j)^2\big)} \; <\infty.
 $$
 The claim now follows directly  from Lemma \ref{lem:order} below.
\endproof

\begin{lemma}\label{lem:order}
Let $r\geq 1$ be an integer and let  $b_1, \ldots, b_r>0$. Assume that $p(x_1,\ldots x_r)\geq 0$ is a non-negative second order polynomial of $r$ real variables such that 
\begin{equation}\label{eq:int}
I(p):=\int_{\R^{r}}\frac{p(x_1,\ldots, x_{r})dx_1\ldots dx_{r} }
{\big(1+x_1^2\big)\ldots \big(1+x_{r}^2\big)\big(1+(\sum_{j=1}^{r}b_jx_j)^2\big)} \; <\infty.
\end{equation}
Then the order of $p$ is at most 2. Conversely, for any polynomial $p$ of second or lower order the above integral is finite.
\end{lemma}

\proof
By performing the change of variables $x'_j:=b_jx_j$ and noting that $(1+x_j^2)\approx(1+(x'_j)^2)$ we may assume that $b_1=\ldots =b_r=1$ .
If $p$ is of second (or lower) order, we have a bound of the form $|p(\vec{x})|\leq C(1+x_1^2+\ldots x_r^2)$. By Fubini's theorem, one checks immediately that $I(1)<\infty$ and $I(x_1^2)<\infty$, which together with symmetry verifies the finiteness of the integral for at most second order polynomials.

Conversely, we assume that $p$ is a non-negative polynomial with degree at least 3. We will show that $I(p)=\infty$. Note that one may assume that $p$ is symmetric, since otherwise it is possible to replace $p$ by $\sum_{\sigma}p(\sigma (x))$, where the sum is over all permutations of the variables $x_1$ to $x_r$. By positivity, this does not decrease the order of $p$.

The case $r=1$ is trivial since then $I(p)=\int_\R p(x_1)(1+x_1^2)^{-2}\diff x_1 = \infty$. Assume then that $r=2$.   If the order of $p$ with respect to one of the variables, say $x_1$, is at least 3, then $x_1\mapsto p(x_1, x_2)$ has order at least 3 for almost every $x_2$, whence $I(p)=\infty$ by Fubini's theorem. Thus $p$ has to be of order 2 separately with respect to both variables, by symmetry and since total  order  is at least 3. We may write 
$$
p(x_1,x_2)=ax_1^2x_2^2+b(x_1+x_2)x_1x_2+ O(|\vec{x}|^2)
$$
Clearly $a$ cannot be strictly negative, and we actually must have $a>0,$ since otherwise $b\not=0$ and we see that $p$ does not stay positive when $\vec{x}$ approaches infinity in one of the directions $\pm (1,1)$.  It follows that on the closed square 
$$
Q_k:= 2^k(1,1)+ 2^{k-3}[-1,1]^2.
$$
with center $(2^k, 2^k)$
we have the estimate $p(x)\gtrsim 2^{4k}$, and the denominator of the integral is bounded by $\lesssim 2^{6k}$. Hence, the integral \eqref{eq:int} over cube $Q_k$ takes values at least of the order $\gtrsim  2^{4k} 2^{-6k}2^{2k} \gtrsim 1$. Since the cubes $Q_k$ are disjoint for $k=1,2,\ldots$, it follows that $I(p)=\infty$.

We take care of the remaining dimensions $r\geq 3$ by applying induction on $r$. Assume that the result in true in dimension $r-1\geq 2$, and consider the case where $p$ is a symmetric and non-negative polynomial on $\R^r$, whose order is at least 3. As before, we see that the maximal order of $p$ with respect to each variable separately is 2, and by denoting $\vec{x'}=:(x_1,\ldots x_{r-1})$ we may write 
\begin{equation}\label{eq:form}
p(x_1,\ldots, x_r)= x_r^2q_2(\vec{x'})+ x_rq_1(\vec{x'}) + q_0(\vec{x'}).
\end{equation}
By Fubini's theorem and noting that for every fixed $x_r$ we have $\big(1+(\sum_{j=1}^{r}b_jx_j)^2\big)\approx \big(1+(\sum_{j=1}^{r-1}b_jx_j)^2\big)$, we see that the assumption that $I(p)<\infty$ leads to
$$
I'(x_r):=\int_{\R^{r-1}}\frac{\big(  x_r^2q_2(\vec{x'})+ x_rq_1(\vec{x'}) + q_0(\vec{x'})\big)dx_1\ldots dx_{r-1} }
{\big(1+x_1^2\big)\ldots \big(1+x_{r-1}^2\big)\big(1+(\sum_{j=1}^{r-1}b_jx_j)^2\big)} <\infty
$$
for almost every $x_r$. At this stage, the inductive assumption yields that $x_r^2q_2(\vec{x'})+ x_rq_1(\vec{x'}) + q_0(\vec{x'})$ is at most a second order polynomial in $\vec{x'}$ for almost every $x_k\in\R.$ This easily implies that the same holds for the individual polynomials $q_j(\vec{x'})$ for $j=0,1,2.$ 

Finally, let us analize what the above conclusion on the polynomials $q_j$ implies in view of representation \eqref{eq:form} and the symmetry of polynomial $p$. By letting $x_r\to\infty$ we see that $q_2\geq 0$. If $q_2$ is not constant, it  has to be a positive symmetric second order polynomial in $\vec{x'}$, and, hence, contains terms $x_j^2$ with positive coefficients. However, then $p$ contains the term $x_1^2x_r^2$ with a positive constant, but not the term $x_1^2x_2^2$. This contradicts the symmetry assumption, and hence $q_2$ is a constant. Then, in order the degree of $p$ to be at least 3, the polynomial $q_2$ has to be of second order, so that the degree of $p$ is 3. But then it cannot be positive, and we have reached the desired contradiction.
\endproof

\begin{coro}
\label{coro:measure_on_lines}
The Nevanlinna  measure $\mu$ on  $\R^n$ in Theorem \ref{thm:hyperplanes} is extremal if and only if  the density $p$ is of the form $p(\vec{x})=(\vec{b}\cdot \vec{x}-c)^2$ 
\end{coro}

\proof
Let $\ell$  be the number of non-zero coefficients $a_k$ in Theorem \ref{thm:hyperplanes}. By choosing suitable coordinates $y_1,\ldots, y_{\ell-1}$ on $M$, it is enough to understand which positive second order polynomials in these are extremals in the cone of positive polynomials. Let $p(y_1,\ldots, y_{\ell-1})$ be a positive polynomial. It may be written in form
$$
p(\vec{y}) = \vec{y}\cdot \mat{A}\vec{y}+ \vec{b} \cdot \vec{y}+c,
$$
where $\mat{A}$ is a positive definite $(\ell-1)\times(\ell-1)$ matrix. We may assume that $\mat{A}$ is diagonal by a rotation. By a dilation of the coordinates we may also assume that the coefficient of each $y_j^2$ is either one or zero, and in case it is zero $p$ does not depend on the $j$-th variable by positivity. Hence, after a translation, we may assume that $p(\vec{y})=c'+\sum_{j=1}^{\ell-1} b_jy_j^2$, where $b_j\in\{0,1\}$. It follows that $\mu$ is not extremal if more that one of the constants $c',b_1,\ldots ,b_\ell$ is zero. This proves one direction. 

In remains to prove the other direction, \ie that constants and squares of first degree polynomials are extremals. It is clear, again after a chance of coordinates, that it is enough to show that both the constant $c$ and the monomial $y_1^2$ are extremals. The case of the constant is evident. Assume then that $y_1^2= p_1(\vec{y})+p_2(\vec{y})$, where $p_k\geq0$,
$k=1,2$. If $p_1$ depends on variable $y_j$ with $j\geq 2$, we obtain a contradiction by letting $y_j\to\infty$. Thus, both $p_1$ and $p_2$ depend only on $y_1$ and the one-dimensional claim is evident.
\endproof

We now give two examples of measures in $\R^n$ where we may discuss their extremality with the help of Theorem \ref{thm:hyperplanes} and Corollary \ref{coro:measure_on_lines}.

\begin{example}
Let $\til{q}(z) = -(z+\alpha)^{-1}$ where $\alpha \in \R$. The representing measure of this function is $\pi\,\delta_\alpha$, \ie the Dirac measure at $\alpha$. This measure is a known extremal measure in dimension one. Let now $k_1,\ldots,k_n \in (0,1]$ be such that $\sum_{j=1}^n k_j = 1$. Then, by \cite[Thm. 4.2]{Nedic2019} and Corollary \ref{coro:measure_on_lines}, the representing measure of the function $q(\vec{z}) = \til{q}(k_1z_1 + \ldots + k_nz_n)$ is extremal on $\R^n$ and has a constant density with respect to the $(n-1)$-dimensional Hausdorff measure on the hyperplane $k_1x_1+\ldots + k_nx_n+\alpha=0$.\hfill$\lozenge$
\end{example}

\begin{example}
Consider the measure $\mu$ on $\R^3$ defined for any Borel set $U \subseteq \R^3$ as
$$\mu(U) := \pi\int_{\R^2}\chi_U(t_1,t_2,-t_1-2t_2)\cdot(2t_1^2+2t_1t_2+t_2^2)\diff t_1 \diff t_2.$$
Clearly, this measure satisfies the growth condition \eqref{eq:measure_growth} and one can verify that this measure is, in fact, a Nevanlinna measure by computing the three intregrals constituting condition \eqref{eq:measure_Nevan}, \eg
$$\int_{\R^2}\frac{1}{(t_1-z_1)^2}\cdot\left(\frac{1}{t_2-z_2} - \frac{1}{t_2 - \overline{z_2}}\right)\cdot\frac{1}{(-t_1-2t_2-\overline{z_3})^2}\diff t_1\diff t_2$$
and two similar others. The \HN function $q$ given by the data $(0,\vec{0},\mu)$ in the sense of Theorem \ref{thm:intRep_Nvar} can then be calculated to be
$$q(z_1,z_2,z_3) = \frac{5 z_1 z_2 + 4 z_1 z_3 + z_2 z_3}{2 (z_1 + 2 z_2 + z_3)}.$$

On the other, alternatively, we  see directly by Theorem \ref{thm:hyperplanes} that this measure in a Nevanlinna measure. Namely, it is supported on the hyperplane $M = \{\vec{t} \in \R^3~|~t_1+2t_2+t_3 = 0\}$ with density $p|_M(\vec{x}) = 2x_1^2+2x_1x_2+x_2^2$ with respect to the orthogonal projection to $t_1$-$t_2$-plane. Observe further that this density may be written, for example, as $p|_M(\vec{x}) = p_1(\vec{x}) + p_2(\vec{x})$ where $p_1(\vec{x}) := x_1^2$ and $p_2(\vec{x}) := (x_1+x_2)^2$, thereby decomposing the measure $\mu$ into two parts $\mu_1$ and $\mu_2$, where 
$$\mu_j(U) := \pi\int_{\R^2}\chi_U(t_1,t_2,-t_1-2t_2)p_j(t_1,t_2)\diff t_1 \diff t_2$$
for $j=1,2$. As before, one can assure themselves by a direct validation of condition \eqref{eq:measure_Nevan} that these measures are indeed Nevanlinna  measures and we note by Corollary \ref{coro:measure_on_lines} that these are extremal measures. The \HN functions represented by the two measure $\mu_1$ and $\mu_2$ can be calculated to be
$$q_1(z_1,z_2,z_3) = \frac{z_1 (2 z_2 + z_3)}{z_1 + 2 z_2 + z_3} \quad\text{and}\quad q_2(z_1,z_2,z_3) = \frac{z_1z_2 + z_2 z_3 + 2 z_1 z_3}{2 (z_1 + 2 z_2 + z_3)}.$$

Finally, note that this is not the only way to decompose the density $p|_M$, as it likewise holds that $p|_M(\vec{x}) = \widetilde{p_1}(\vec{x}) + \widetilde{p_2}(\vec{x})$ where $\widetilde{p_1}(\vec{x}) := (\sqrt{2}x_1+\tfrac{1}{\sqrt{2}}x_2)^2$ and $\widetilde{p_2}(\vec{x}) := \tfrac{1}{2} x_2^2$. Using these densities to define two measure $\widetilde{\mu_1}$ and $\widetilde{\mu_2}$ in an analogous way as before, one arrives at a different decomposition of the measure $\mu$ as a sum of extremal measures.\hfill$\lozenge$
\end{example}

We conclude this section with the following observation that strengthens the conclusion of Theorem \ref{thm:hyperplanes} for hyperplanes that cannot carry a non-trivial Nevanlinna measure: in fact the restriction of any  Nevanlinna measure  onto such a  hyperplane vanishes.

\begin{thm}
\label{thm:lines_with_positive_slope} 
Assume that $L\subset \R^n$ is a $(n-1)$-dimensional hyperplane whose normal is not contained in $C_-\cup C_+$ (see formula \eqref{eq:perus}). Then $\mu(L)=0$ for all Nevanlinna measures on $\R^n$. 
\end{thm}

\proof
Assume that $\mu$ and $L$ are as in the statement. By translation, we may assume that $L$ contains the origin and if $\vec{a}\in\R^n$ is a normal vector of $L$, then $\vec{a}$ is not contained in $C_-\cup C_+$. Choose orthogonal coordinates $\vec{y}$ so that $y_n$-axis has  direction $\vec{a}$. Denote $(y_1,\ldots , y_{n-1})=\vec{y}'$ so that $\vec{y}=(\vec{y}',y_n)$. Pick a smooth one-dimensional test function $\psi_0$ supported on $(-1,1)$ such that $\widehat \psi_0>0$ everywhere on $\R$ with $\widehat\psi(0)=1$, and write 
$$
\psi_1(\vec{y}'):=\psi_0(y_1)\ldots \psi_0 (y_{n-1}) \quad \textrm{with}\quad \widehat \psi_1(\vec{\xi}')= \widehat \psi_0(\xi_1)\ldots \widehat \psi_0 (\xi_{n-1}),
$$
where we denoted by $\vec{\xi}=(\vec{\xi}',\xi_n)$ the Fourier variable corresponding to coordinates $\vec{y}$.
It is enough to verify that
$$
\int_L\widehat \psi_1(\vec{y}')\diff\mu_{|L}(\vec{y}')=0.
$$
 It follows that
\begin{eqnarray*}
\int_L\widehat \psi_1(\vec{y}')\diff\mu_{|L}(\vec{y}')&=&\lim_{\varepsilon\to 0}\int_{\R^n}\widehat\psi_1(\vec{y}')\widehat\psi_0(y_n/\varepsilon)\diff\mu(\vec{y})\\
&=&\lim_{\varepsilon\to 0}\langle\varepsilon \psi_1(\vec{\xi}')\psi_0(\varepsilon \xi_n)\, ,\, \widehat \mu\rangle \\
&=&0.
\end{eqnarray*}
The final conclusion above was obtained by observing first that since $\psi_1$ has compact support in $\vec{y'}$ and the $y_n$-axis intersects $C_-\cup C_+$ only at the origin, we have for some $R_0>0$ such that
$$
\supp\big(\varepsilon \psi_1(\vec{\xi}')\psi(\varepsilon \xi_n)\big)\; \cap \supp(\widehat\mu)\; \subset \; B(\vec{0},R_0).
$$
Moreover, we note that 
$\varepsilon \psi_1(\vec{y}')\psi(\varepsilon y_n)\to 0$   together with all of its derivatives uniformly (say) in the ball $B(\vec{0},2R_0)$.
This finishes the proof.
\endproof

\section{Dependence of the parameters on the fixed variables}
\label{sec:fixed_var}

A well known way of constructing new \HN functions out of old ones is to \eg replace an independent variable with a non-negative linear combination of independent variables. One case where such a method is particularly fruitful arises when replacing the independent variable of a \HN function in $\C^+$ with a convex combination of arbitrary many independent variables \cite{Nedic2019}. In particular, knowing the representing parameters of the original function one can directly explicitly write down the representing parameters of the new function \cite[Thm. 4.2]{Nedic2019}. In what follows, we will consider a problem going in the converse direction, namely to understand how the representing parameters of \HN function change when some of the independent variables are set to fixed values.

\subsection{Introduction to the problem}

Let $n,m \in \N$ and let $q$ be a \HN function of $n+m$ variables. From this function we construct a new function of only $n$ variables by fixing the remaining $m$ variables. More precisely, let $\vec{\zeta} \in \C^{+m}$ be fixed and let $\til{q}_{\vec{\zeta}}$ be the the \HN function of $n$ variables defined as
$$\til{q}_{\vec{\zeta}}\colon \vec{z} \mapsto q(\vec{z},\vec{\zeta}).$$
Let $(\til{a}_{\vec{\zeta}},\tilvec{b}_{\vec{\zeta}},\til{\mu}_{\vec{\zeta}})$ be the representing parameters of the function $\til{q}_{\vec{\zeta}}$ in the sense of Theorem \ref{thm:intRep_Nvar}. Here, the dependence of the measure $\til{\mu}_{\vec{\zeta}}$ on the value of $\vec{\zeta} \in \C^{+m}$ is of particular interest as it describes the limit of the function $q$ at a particular part of the non-distinguished boundary of the poly-upper half-plane.

We begin by presenting three examples that highlight some particularly interesting cases of the dependence of the measure $\til{\mu}_{\vec{\zeta}}$ on $\vec{\zeta} \in \C^{+m}$.

\begin{example}\label{ex:fixed_var3}
Let $q(z_1,z_2) := -(z_1+z_2+\I)^{-1}$. The representing parameters $(a,(b_1,b_2),\mu)$ of this function in the sense of Theorem \ref{thm:intRep_Nvar} can be shown to be
$$a = b_1 = b_2 = 0,$$
while the measure $\mu$ is defined for $U \subseteq \R^2$ as
$$\mu(U) = \int_{\R^2}\frac{\chi_U(t_1,t_2)}{1+t_1^2+t_2^2}\diff t_1 \diff t_2,$$
which is absolutely continuous with respect to $\lambda_{\R^2}$.

Let now $\til{q}_{\zeta}(z) := q(z,\zeta)$ for a fixed $\zeta \in \C^+$. Then, it holds that
$$\til{q}_{\zeta}(z) = -\frac{1}{z+\zeta + \I}$$
with
$$\til{a}_\zeta = \Re[\til{q}_\zeta(\I)] = - \frac{\Re[\zeta]}{|2\,\I + \zeta|^2} \quad\text{and}\quad \til{b}_\zeta = \lim\limits_{z \ntto \infty}\frac{\til{q}_\zeta(z)}{z} = 0.$$
Furthermore, it holds that
$$\Im[\til{q}_\zeta(x+\I\,y)] = \frac{y+\Im[\zeta]+1}{(x+\Re[\zeta])^2 + (y + \Im[\zeta]+1)^2}.$$
Hence, the measure $\til{\mu}_\zeta$ is absolutely continuous with respect to $\lambda_\R$ with its density being
$$t \mapsto \frac{\Im[\zeta]+1}{(t+\Re[\zeta])^2 + (\Im[\zeta]+1)^2} = \Im[-(t+\zeta+\I)^{-1}].$$
Note that this density is a harmonic function in the parameter $\zeta$, as are the functions $\zeta \mapsto \til{a}_\zeta$ and $\zeta \mapsto \til{b}_\zeta$. \hfill$\lozenge$
\end{example}

\begin{example}\label{ex:fixed_var1}
Let $q(z_1,z_2) := -(z_1+z_2)^{-1}$. Similarly to the function $q$ from the previous example, it holds for this function that $a = b_1 = b_2 = 0$. However, the representing measure $\mu$ of this function is quite different, namely it is defined for $U \subseteq \R^2$ as
$$\mu(U) = \pi\int_\R\chi_U(t,-t)\diff t,$$
which is singular continuous with respect to $\lambda_{\R^2}$.

Despite the major difference in the properties of the representing measure compared to the previous example, the properties of the functions $\til{q}_{\zeta}$ and from this and the previous example are quite similar. Indeed, let $\til{q}_{\zeta}(z) := q(z,\zeta)$ for a fixed $\zeta \in \C^+$ as before. Then, it holds that
$$\til{q}_{\zeta}(z) = -\frac{1}{z+\zeta} \quad\text{with}\quad \til{a}_\zeta = - \frac{\Re[\zeta]}{|\I + \zeta|^2} \quad\text{and}\quad \til{b}_\zeta = 0.$$
Furthermore, it holds that
$$\Im[\til{q}_\zeta(x+\I\,y)] = \frac{y+\Im[\zeta]}{(x+\Re[\zeta])^2 + (y + \Im[\zeta])^2}.$$
Hence, as in Example \ref{ex:fixed_var3}, the measure $\til{\mu}_\zeta$ is absolutely continuous with respect to $\lambda_\R$ with its density being
$$t \mapsto \frac{\Im[\zeta]}{(t+\Re[\zeta])^2 + \Im[\zeta]^2} = \Im[-(t+\zeta)^{-1}].$$
Additionally, as previously, the density is a harmonic function in the parameter $\zeta$, as are the functions $\zeta \mapsto \til{a}_\zeta$ and $\zeta \mapsto \til{b}_\zeta$. \hfill$\lozenge$
\end{example}

\begin{example}\label{ex:fixed_var2}
Let $q(z_1,z_2) := -z_1^{-1} - z_2^{-1}$. Also for this function $q$, as the previous two, it holds that $a = b_1 = b_2 = 0$. The representing measure $\mu$ this time equals
and 
$$\mu = \pi\delta_0 \otimes \lambda_\R + \lambda_\R \otimes \pi\delta_0,$$
which, like the measure in Example \ref{ex:fixed_var1}, is singular continuous with respect to $\lambda_{\R^2}$, though this time it is a particularly a sum of two product measures. 

Let again $\til{q}_\zeta(z) := q(z,\zeta)$ for a fixed $\zeta \in \C^+$. Then, it holds that
$$\til{q}_\zeta(z) = -\frac{1}{z} - \frac{1}{\zeta} \quad\text{with}\quad \til{a}_\zeta = - \frac{\Re[\zeta]}{|\zeta|^2} \quad\text{and}\quad \til{b}_\zeta = 0.$$
Furthermore, it holds that
$$\Im[\til{q}_\zeta(x+\I\,y)] = \frac{y}{x^2+y^2} + \frac{\Im[\zeta]}{\Re[\zeta]^2 + \Im[\zeta]^2}.$$
Hence, the measure $\til{\mu}_\zeta$ equals
$$\til{\mu}_\zeta = \pi\delta_0 + \frac{\Im[\zeta]}{\Re[\zeta]^2 + \Im[\zeta]^2}\lambda_\R = \pi\delta_0 + \Im[-\zeta^{-1}]\lambda_\R.$$
Note that, here, the measure $\til{\mu}_\zeta$ can be written as a sum of a measure depending on $\zeta$ and a measure that is independent of $\zeta$. The one depending on $\zeta$ is absolutely continuous with respect to $\lambda_\R$ with a density that is harmonic in the parameter $\zeta$ while the independent part is singular with respect to $\lambda_\R$. Again, the functions $\zeta \mapsto \til{a}_\zeta$ and $\zeta \mapsto \til{b}_\zeta$ are also harmonic in the variable $\zeta$.\hfill$\lozenge$
\end{example}

The above examples show also another interesting phenomenon. If one wishes to consider the case when the fixed variable $\zeta$ is taken at the boundary, then nothing can be said in general. For the function form Example \ref{ex:fixed_var3}, the corresponding function $\til{q}_\zeta$ has a representing measure that is absolutely continuous even for any fixed $\zeta \in \R$. On the other hand, for the function form Example \ref{ex:fixed_var1}, the measure of the corresponding function $\til{q}_\zeta$ becomes a pure point measure if we fix $\zeta \in \R$. Moreover, by considering the function $(z_1,z_2) \mapsto -(z_1 + \sqrt{z_2}))^{-1}$, we see that when the other variable is fixed from the boundary, the representing measure with respect to the variable $z_1$ can be singular for some $z_2\in\R$ and absolutely continuous for other $z_2\in\R$. Additionally, for the function from Example \ref{ex:fixed_var2}, the corresponding function $\til{q}_\zeta$ is not well-defined for $\zeta = 0$.  

Using the results from the previous section, we can, however, describe 
the properties of the measure $\til{\mu}_\zeta$ of the function $\til{q}_\zeta$ for functions whose measure is supported on a line with negative slope in $\R^2$. Indeed, let $\mu$ be as stated and let $q$ be the \HN function given by the data $(0,\vec{0},\mu)$ in the sense of Theorem \ref{thm:intRep_Nvar}. Without loss of generality, we may restrict ourselves to the case when $\mu$ is supported on a line with negative slope through the origin. All other cases may be handled via translations. Furthermore, we may always assume that we are setting the second variable to a fixed value. By Corollary \ref{coro:measure_on_lines}, we know exactly what all such measures $\mu$ are, implying that the function $q$ is of the form
$$q(z_1,z_2) = \frac{\gamma z_1 z_2 - \eta}{\alpha z_1 + \beta z_2}$$
for some constants $\alpha,\beta > 0$ and $\gamma,\eta \geq 0$. It now follows by standard one-variable theory that for every point $\zeta \in \C^+$ the measure $\til{\mu}_\zeta$ is absolutely continuous with respect to $\lambda_\R$, while for every point $\zeta \in \R$ the measure $\til{\mu}_\zeta$ is well-defined and is a pure point measure.

In the next section, we will focus on results that are valid in greater generality than the above special case in dimension two. In particular, Theorems \ref{thm:singular_part_polydisc} and \ref{thm:singular_part_HN} show that what occurs in Example \ref{ex:fixed_var2} can be considered a generic case. Additionally, Corollary \ref{coro:absolutely_continuous} show that what occurs in Example \ref{ex:fixed_var3} holds for all \HN functions with absolutely continuous measures.

\subsection{Dependence of the singular part on the fixed parameter}

We start by considering first the situation in the case of the polydisc. To begin with, we first consider measures on $\T^n$ that are parametrized by the $\vec{\zeta} \in \D^m$ in such a way that the dependence is harmonic.
We denote by $M(\T^n)$ he Banach space of of signed and finite Borel measures on $\T^n$.

\begin{lemma}\label{lem:measure_polydisc}
Let $f:\D^m \to M(\T^n)$ be a separately harmonic (resp. pluriharmonic) and bounded map. In other words,  we assume that 
$$
\sup_{\vec{\zeta}\in\D^m} \|f(\vec{\zeta})\|_{M(\T^n)}<\infty,
$$
and separate harmonicity (resp. pluriharmonicity) is understood in the sense that for all $\varphi\in C^\infty (\T^n)$ the function 
$$
\vec{\zeta} \mapsto \int_{\T^n} \varphi(\vec{w})\diff f(\vec{\zeta})(\vec{w})
$$
is separately harmonic (resp. pluriharmonic) on $\D^m$. Then 
\begin{itemize}
    \item[(i)]{The function $\vec{\zeta}\mapsto f(\vec{\zeta})(E)$ is separately harmonic (resp. pluriharmonic) on $\D^m$ for all Borel-sets $E\subset\T^m$.}
    \item[(ii)]{There is a Borel set $A\subset \T^n$ of measure zero, a bounded positive singular measure $\nu$ supported on $A$, and pointwise separately harmonic (resp. pluriharmonic) functions (Borel in the space variable) $h_{\vec{\zeta}}\colon \T^n \to \C$ and $g_{\vec{\zeta}}\colon A \to \C$ with $\|h_{\vec{\zeta}}\|_{L^1(\T^n)}\leq C$ and $\|g_{\vec{\zeta}}\|_{L^1(\nu)}\leq C$ for all ${\vec{\zeta}}\in\D^m$, and such that
    $$\diff f(\vec{\zeta})(\vec{w}) = h_{\vec{\zeta}}(\vec{w})\diff\vec{w} + g_{\vec{\zeta}}(\vec{w})\diff\nu(\vec{w}) \qquad \textrm{for all}\quad \vec{\zeta} \in \D^m.$$
    }
\end{itemize}
\end{lemma}
\proof
From the assumption it follows that  $f\colon \D^m\to M(\T^n)$ is separately harmonic (resp. pluriharmonic) in the classical sense \cite{Arendt2016} as a Banach space valued function, and it has a Poisson integral representation on each polydisc $\overline {B(0,r)}^m$. The pluriharmonicity is naturally defined again by demanding that $f$ is a harmonic function on each complex line. Especially, since evaluation on a Borel set belongs to the dual of $M(\T^n)$,  statement (i) follows. In addition, we obtain\footnote{This is obtained just as the proof of the power series expansion in  \cite[Thm. 5.2]{Arendt2016} using analogue of \cite[Formula (5)]{Arendt2016} where the Poisson formula is replaced by the Poisson-formula in the polydisc of radius $r<1$, and the estimate for the Fourier coefficients is obtained by letting $r\to 1^-$.}  that for all multi-indices $\vec{k}\in \Z^m$ there are uniformly  bounded measures $\mu_{\vec{k}}\in M(\T^n)$ such that
\begin{equation}
\label{eq:measure_series_decomposition}
f({\vec{\zeta}}) = \sum_{\vec{k}\in\Z^m}\mu_{\vec{k}}\vec{\zeta}^{\vec{k}},
\end{equation}
where $\vec{\zeta}^{\vec{k}} := (\zeta_1)^{k_1} \ldots (\zeta_m)^{k_m}$ and we interpret $(\zeta_j)^\ell= (\overline{\zeta_j})^{|\ell|}$ for negative integers $\ell$. Decompose now
$$\diff\mu_{\vec{k}}(\vec{w}) = H_{\vec{k}}(\vec{w})\diff\vec{w} + \diff\nu_{\vec{k}}(\vec{w})$$
into its absolutely continuous and singular parts, and write $\nu := \sum_{\vec{k}\in\Z^m}2^{-|\vec{k}|}\nu_{\vec{k}}$  so that $\nu$ is singular. Write $\diff\nu_{\vec{k}}(\vec{w})=G_{\vec{k}}(\vec{w})\diff\nu(\vec{w})$. One easily checks  that the functions
\begin{equation}\label{eq:-1}
h_{\vec{\zeta}}(\vec{w}) :=  \sum_{\vec{k} \in \Z^m} \vec{\zeta}^{\vec{k}}H_{\vec{k}}(\vec{w}) \qquad\textrm{and}\qquad
g_{\vec{\zeta}}(\vec{w}) :=  \sum_{\vec{k} \in \Z^m} \vec{\zeta}^{\vec{k}}G_{\vec{k}}(\vec{w})
\end{equation}
do the job, because  the geometric convergence in $\zeta$ implies pointwise absolute convergence in the definition of $h_{\vec{\zeta}}$ at all points $\vec{w}$ outside the nullset 
$$
\bigcap_{j=1}^\infty\bigcup_{|\vec{k}|\geq j}\big\{\vec{w}\in\T^n\,\big|\, |H_{\vec{k}}(\vec{w})|\geq |k|^{m+1}\big\},
$$
and a similar argument  works the singular part. This implies part (ii) in the separately harmonic case. Finally, in the pluriharmonic case, the fact that the function $\vec{\zeta}\mapsto \int_{\T^n} \varphi(\vec{w})\diff f(\vec{\zeta})(\vec{w})$ is pluriharmonic for all $\varphi\in C^\infty (\T^n)$ implies that  $\int_{\T^n} \varphi(\vec{w})\diff \mu_{\vec{k}}(\vec{w})$ vanishes  if some components of $k$ have different sign. Thus $\mu_k$ vanishes for such $k$, and this clearly implies the plurisubharmonicity of functions $h_{\vec{\zeta}}$ and $g_{\vec{\zeta}}$ via \eqref{eq:measure_series_decomposition} and \eqref{eq:-1}.
\endproof

\begin{prop}\label{prop:singular_part_polydisc}
Assume that $f$ is as in the previous lemma but with values on positive measures. Assume that one of the measures, say  $\mu_{\vec{\zeta}_0}$, has a non-trivial singular part. Then $f(\vec{\zeta})$ has a non-trivial singular part for all $\vec{\zeta}\in\D^m$.
\end{prop}

\proof
Let $A\subset\T^n$ be as in the previous proof. Simply apply 1-dimensional Harnack inequality (see \cite[Sec. 6.3.2]{Ahlfors1979}),  or alternatively the local minimum principle, recursively in each variable on the non-negative separately harmonic function $\vec{\zeta} \mapsto \mu_{\vec{\zeta}}(A)$.
\endproof

One should note that above the measures were not assumed to be RP-measures. For them, we have a stronger result, which was our main goal in this subsection. We state it in terms of RP-functions as part (ii) of the following result.

\begin{thm}
\label{thm:singular_part_polydisc}
Let $G$ be a RP-function on $\D^{n+m}$. Denote by $\nu$ the representing RP-measure of $G$ in the sense of Theorem \ref{thm:RP_intRep_Nvar}. For $\vec{\zeta} \in \D^m$ denote by $\nu_{\vec{\zeta}}$ the representing measure of the RP-function $G_{\vec{\zeta}}\colon\D^n\to \C$, where
$$
G_{\vec{\zeta}}(z_1,\ldots, z_n):= G(z_1,\ldots, z_n,\zeta_1,\ldots,\zeta_m).
$$ 
{\rm (i)} \quad The map $\zeta \mapsto \nu_\zeta$ is plurisubharmonic in the sense of {\rm Lemma  \ref {lem:measure_polydisc}}, and in particular $\zeta \mapsto \nu_\zeta(A)$ is pluriharmonic in $\D^m$ for any Borel set $A\subset \T^n.$

\smallskip

\noindent {\rm (ii)}\quad  if  $\nu_\zeta$ has a non-trivial singular part for one $\vec{\zeta}\in\D^m$, then it has a non-trivial singular part for all other values and the singular part is independent of $\zeta$. Furthermore, we may write
$$
G(\vec{z},\vec{\zeta}) = G_1(\vec{z}) + G_2(\vec{z},\vec{\zeta})$$
for all $\vec{z}\in\D^n$ and $\vec{\zeta}\in \D^m$, where both $G_1$ and $G_2$ are RP-functions and the representing measure of $G_2(\,\cdot\,,\zeta)$ is absolutely continuous for all $\vec{\zeta} \in \D^m$.
\end{thm}

\proof
(i) Given $\varphi\in C(\T^n)$ we note that $\zeta\to \int_{\T^n}G((1-1/k)z,\zeta)\varphi(z)d\lambda_{\T^n}(z)$ is plurisubharmonic and uniformly bounded in $k,$ for $|\xi|<r<1$, for any $r<1$. Letting $k\to\infty$ we obtain that the limit $\zeta\mapsto \int_{\T^n} \varphi(\vec{w})\diff \nu_\zeta(\vec{w})$ is pluriharmonic as a pointwise limit of uniformly bounded pluriharmonic functions.

(ii) According to part (i) and Lemma \ref {lem:measure_polydisc} we may  decompose the measure $\nu_{\vec{\zeta}}$ as in \eqref{eq:measure_series_decomposition}.  Now, write further each of the measures $\nu_{\vec{k}}$ in terms of their Fourier series, \ie $\nu_{\vec{k}}=\sum_{\vec{j}\in\Z^n} a_{\vec{j},\vec{k}}\E^{\I\,\vec{j}\dot w},\,a_{\vec{j},\vec{k}} \in \C,$ where the convergence is in the weak$*$ and substitute it back into decomposition \eqref{eq:measure_series_decomposition} of $\nu_{\vec{\zeta}}$. Using our notational convention for negative powers as in the proof Lemma \ref{lem:measure_polydisc}, we obtain
$$
\Re[G(\vec{z},\vec{\zeta})] = \sum_{\substack{\vec{j}\in\Z^n \\ \vec{k}\in \Z^m}}a_{\vec{j},\vec{k}}\vec{z}^{\vec{j}}\vec{\zeta}^{\vec{k}}.
$$
This implies that $a_{\vec{j},\vec{k}}=0$ if a component of $\vec{j}$ have different sign than some component of $\vec{k}$. Especially, for $\vec{k}\neq\vec{0}$, the measure $\nu_{\vec{k}}$ is either analytic (i.e. all non-zero Fourier coefficients have only non-negative indices) or anti-analytic (i.e. all non-zero Fourier coefficients have only non-positive indices).
In particular, it is absolutely continuous by the brothers Riesz theorem \cite[Thm. 17.13]{Rudin1987} and, hence, $\nu_{\vec{0}}$ is the only measure which may be singular. The required decomposition is readily obtained by constructing $G_1$ from the singular part of $\mu_0$.
\endproof

By applying the correspondence between \HN functions and RP-functions as recalled in Section \ref{subsec:correspondence}, in particular using formula \eqref{eq:nu_to_mu}, we may write a version of Theorem \ref{thm:singular_part_polydisc} for \HN functions.

\begin{thm}\label{thm:singular_part_HN}
Let $q$ be a \HN function in $\C^{+(n+m)}$. For $\zeta\in \C^{+m}$ denote let $\til{\mu}_{\vec{\zeta}}$ denote the representing measure of the \HN function $\til{q}_{\vec{\zeta}}\colon\C^{+n}\to \C$, where
$$
\til{q}_{\vec{\zeta}}(z_1,\ldots, z_n) := q(z_1,\ldots, z_n,\zeta_1,\ldots,\zeta_m).
$$ 
Then, if $\til{\mu}_{\vec{\zeta}}$ has a non-trivial singular part for one $\vec{\zeta}\in\C^{+m}$, it has a non-trivial singular part for all other values as well, and the singular part is independent of $\zeta$. Furthermore, we may write
$$
q(\vec{z},\vec{\zeta}) = q_1(\vec{z}) + q_2(\vec{z},\vec{\zeta})
$$
for all $\vec{z} \in \C^{+n}$ and $\vec{\zeta} \in \C^{+m}$, where both $q_1$ and $q_2$ are \HN functions and the representing measure of $q_2(\,\cdot\,,\vec{\zeta})$ is absolutely continuous for all $\vec{\zeta}\in\C^{+m}$.
\end{thm}

\begin{remark}
The statements of Theorems \ref{thm:singular_part_polydisc} and \ref{thm:singular_part_HN} may easily be adapted to the case when the independent and fixed variables are arbitrarily permuted.
\end{remark}

Theorems \ref{thm:singular_part_polydisc} and \ref{thm:singular_part_HN} immediately imply the following corollary which we only state in the poly-upper half-plane as explicit examples of absolutely continuous measures are more often presented in $\R^n$, see \eg Example \ref{ex:Fourier_examples_dim2}.

\begin{coro}
\label{coro:absolutely_continuous}
Let $\mu$ be a Nevanlinna measure that is absolutely continuous with respect to $\lambda_{\R^{n+m}}$. Then, for every point $\vec{\zeta} \in \C^{+m}$, the measure $\til{\mu}_{\vec{\zeta}}$ is absolutely continuous with respect to $\lambda_{\R^n}$.
\end{coro}

\proof
If $\til{\mu}_{\vec{\zeta}}$ would have a non-trivial singular part for some $\vec{\zeta} \in \C^{+m}$, the starting measure $\mu$ would had to have had a non-trivial singular part also, contradicting our starting assumption.
\endproof

\subsection{Refinements using the Stieltjes inversion formula}

One of the major advantages one has at their disposal when working in the poly-upper half-plane as opposed to the unit polydisc is the elegance and simplicity with which one can describe the representing measure of \HN functions using the Stieltjes inversion formula \eqref{eq:Stieltjes_inversion}. With its help, the following proposition describes in general how the representing parameters $(\til{a}_{\vec{\zeta}},\tilvec{b}_{\vec{\zeta}},\til{\mu}_{\vec{\zeta}})$ of the function $\til{q}_{\vec{\zeta}}$ depend on the fixed variables $\vec{\zeta}\in \C^{+m}$.

\begin{prop}
\label{prop:tilde_parameters}
Let $q$ be a \HN function in $n+m$ variables with $n,m \in \N$. Let $\vec{\zeta} \in \C^{+m}$ be fixed and let $\til{q}_{\vec{\zeta}}$ be the the \HN function of $n$ variables defined as
$$\til{q}_{\vec{\zeta}}\colon \vec{z} \mapsto q(\vec{z},\vec{\zeta}).$$
Then, it holds for the representing parameters $(\til{a}_{\vec{\zeta}},\tilvec{b}_{\vec{\zeta}},\til{\mu}_{\vec{\zeta}})$ of the function $\til{q}_{\vec{\zeta}}$ that
\begin{equation}
    \label{eq:a_and_b_tilde}
    \til{a}_{\vec{\zeta}} = a + \sum_{\ell = 1}^mb_{n+\ell}\,\Re[\zeta_\ell], \quad \tilvec{b}_{\vec{\zeta}} = (b_1,\ldots,b_n)
\end{equation}
and for any function $\psi$ as in the Stieltjes inversion formula \eqref{eq:Stieltjes_inversion} it holds that
\begin{multline}
\label{eq:mu_tilde}
    \int_{\R^n}\psi(\vec{t})\diff\til{\mu}_{\vec{\zeta}}(\vec{t}) \\ = \left(\sum_{\ell = 1}^mb_{n+\ell}\,\Im[\zeta_\ell]\right)\int_{\R^n}\psi(\vec{x})\diff \vec{x} + \frac{1}{\pi^m}\int_{\R^{n+m}}\psi(\vec{t})\pois_m(\vec{\zeta},\vec{\tau})\diff\mu(\vec{t},\vec{\tau}).
\end{multline}
\end{prop}

\proof
Before we begin, note that within this proof, it always holds that $\vec{z} = \vec{x} + \I\,\vec{y} \in \C^{+n}$ and $\vec{t} \in \R^n$ while $\vec{\zeta} \in \C^{+m}$ and $\vec{\tau} \in \R^{m}$. Moreover, for $k \in \N$, we denote
\begin{equation}
    \label{eq:1_vec}
    \vec{1}_k := \underbrace{(1,1,\ldots,1)}_{k \text{ entries}} \in \C^k.
\end{equation}

The integral representation of the function $\til{q}_{\vec{\zeta}}$ as a \HN function of $n$ variables implies, for $\vec{z} \in \C^{+n}$, that
\begin{multline*}
    \til{q}_{\vec{\zeta}}(\vec{z}) = \til{a}_{\vec{\zeta}} + \sum_{j = 1}^n(\tilvec{b}_{\vec{\zeta}})_j\,z_j + \frac{1}{\pi^n}\int_{\R^n} K_n(\vec{z},\vec{t})\diff\til{\mu}_{\vec{\zeta}}(\vec{t}) \\ 
    = \til{a}_{\vec{\zeta}} + \sum_{j = 1}^n(\tilvec{b}_{\vec{\zeta}})_j\,z_j + \frac{1}{\pi^{n+m}}\int_{\R^{n+m}} K_{n+m}((\vec{z},\I\,\vec{1}_m),(\vec{t},\vec{\tau}))\diff(\til{\mu}_{\vec{\zeta}} \otimes \lambda_{\R^m})(\vec{t},\vec{\tau}).
\end{multline*}
On the other hand, using the integral representation of the function $q$ to describe the function $\til{q}_{\vec{\zeta}}$, we get that
$$\til{q}_{\vec{\zeta}}(\vec{z}) \\ = a + \sum_{j=1}^nb_j\,z_j + \sum_{\ell = 1}^mb_{n+\ell}\,\zeta_\ell + \frac{1}{\pi^{n+m}}\int_{\R^{n+m}}K_{n+m}((\vec{z},\vec{\zeta}),(\vec{t,\vec{\tau}}))\diff\mu(\vec{t},\vec{\tau}).$$
Observe now that for any $\ell \in \{1,\ldots,m\}$ it holds that
\begin{eqnarray*}
b_{n+\ell}\,\zeta_\ell & = & b_{n+\ell}\,\Re[\zeta_\ell] + \I\,b_{n+\ell}\,\Im[\zeta_\ell] \\[0.2cm]
~ & = & b_{n+\ell}\,\Re[\zeta_\ell] + \frac{b_{n+\ell}\,\Im[\zeta_\ell]}{\pi^{n+m}}\int_{\R^{n+m}}K_{n+m}((\I\,\vec{1}_n,\I\,\vec{1}_m),(\vec{t},\vec{\tau}))\diff \vec{t} \, \diff \vec{\tau}.
\end{eqnarray*}
Comparing the two representations of the function $\til{q}_{\vec{\zeta}}$, we infer, via uniqueness of the representing parameters, that the parameters $\til{a}_{\vec{\zeta}}$ and $\tilvec{b}_{\vec{\zeta}}$ are indeed described by formula \eqref{eq:a_and_b_tilde}. Furthermore, we infer that the measure $\til{\mu}_{\vec{\zeta}}$ is the unique positive Borel measure on $\R$ that solves the equation
\begin{multline*}
    \int_{\R^{n+m}} K_{n+m}((\vec{z},\I\,\vec{1}_m),(\vec{t},\vec{\tau}))\diff(\til{\mu}_{\vec{\zeta}} \otimes \lambda_{\R^m})(\vec{t},\vec{\tau}) \\ = \int_{\R^{n+m}}K_{n+m}((\vec{z},\vec{\zeta}),(\vec{t},\vec{\tau}))\diff\mu(\vec{t},\vec{\tau}) \\ + \left(\sum_{\ell=1}^mb_{n+\ell}\Im[\zeta_\ell]\right)\cdot\int_{\R^{n+m}}K_{n+m}((\I\,\vec{1}_n,\I\,\vec{1}_m),(\vec{t},\vec{\tau}))\diff\lambda_{\R^{n+m}}(\vec{t},\vec{\tau}).
\end{multline*}
Since all of the measures in the above equality are Nevanlinna measures, we may take the imaginary part of the above equality and invoke  the Stieltjes inversion formula \eqref{eq:Stieltjes_inversion} to obtain
$$\begin{array}{RCL}
\multicolumn{3}{L}{\int_{\R^n}\psi(\vec{t})\diff\til{\mu}_{\vec{\zeta}}(\vec{t}) = \lim\limits_{\vec{y} \to \vec{0}^+}\int_{\R^n}\psi(\vec{x})\Im[\til{q}_{\vec{\zeta}}(\vec{x}+\I\,\vec{y})]\diff \vec{x}} \\[0.4cm]
~ & = & \lim\limits_{\vec{y} \to \vec{0}^+}\int_{\R^n}\psi(\vec{x})\bigg(\sum_{j=1}^n(\tilvec{b}_{\vec{\zeta}})_jy_j + \frac{1}{\pi^n}\int_{\R^n} \pois_n(\vec{x}+\I\,\vec{y},\vec{t})\diff\til{\mu}_{\vec{\zeta}}(\vec{t})\bigg)\diff \vec{x} \\[0.4cm]
~ & = & 0 + \lim\limits_{\vec{y} \to \vec{0}^+}\int_{\R^n}\psi(\vec{x})\bigg(\frac{1}{\pi^{n+m}}\int_{\R^{n+m}} \pois_{n+m}((\vec{x}+\I\,\vec{y},\I\,\vec{1}_m),(\vec{t},\vec{\tau}))\diff\til{\mu}_{\vec{\zeta}}(\vec{t})\diff\vec{\tau}\bigg)\diff \vec{x} \\[0.4cm]
~ & = & \lim\limits_{\vec{y} \to \vec{0}^+}\int_{\R^n}\psi(\vec{x})\bigg(\frac{1}{\pi^{n+m}}\int_{\R^{n+m}} \pois_{n+m}((\vec{x}+\I\,\vec{y},\vec{\zeta}),(\vec{t},\vec{\tau}))\diff\mu(\vec{t},\vec{\tau}) \\[0.4cm]
~ & ~ & \quad + \left(\sum_{\ell = 1}^mb_{n+\ell}\,\Im[\zeta_\ell]\right)\frac{1}{\pi^{n+m}}\int_{\R^{n+m}}\pois_{n+m}((\I\,\vec{1}_n,\I\,\vec{1}_m),(\vec{t},\vec{\tau}))\diff\vec{t}\,\diff\vec{\tau}\bigg)\diff \vec{x} \\[0.4cm]
~ & = & \left(\sum_{\ell = 1}^mb_{n+\ell}\,\Im[\zeta_\ell]\right)\int_{\R^n}\psi(\vec{x})\diff \vec{x} \\[0.4cm]
~ & ~ & \quad + \lim\limits_{\vec{y} \to \vec{0}^+}\frac{1}{\pi^{n+m}}\int_{\R^{n+m}}\pois_m(\vec{\zeta},\vec{\tau})\bigg(\int_{\R^n}\pois_n(\vec{x}+\I\,\vec{y},\vec{t})\psi(\vec{x})\diff \vec{x}\bigg)\diff\mu(\vec{t},\vec{\tau})  \\[0.4cm]
~ & = & \left(\sum_{\ell = 1}^mb_{n+\ell}\,\Im[\zeta_\ell]\right)\int_{\R^n}\psi(\vec{x})\diff \vec{x} + \frac{1}{\pi^m}\int_{\R^{n+m}}\psi(\vec{t})\pois_m(\vec{\zeta},\vec{\tau})\diff\mu(\vec{t},\vec{\tau}),
\end{array}$$
as desired.
\endproof

It is clear from the above proposition that the functions $\vec{\zeta} \mapsto \til{a}_{\vec{\zeta}}$ and $\vec{\zeta} \mapsto \tilvec{b}_{\vec{\zeta}}$ are pluriharmonic functions on $\C^{+m}$. The following proposition now establishes an analogous statement for the function $\vec{\zeta} \mapsto \til{\mu}_{\vec{\zeta}}$.

\begin{prop}
Let $q$ be a \HN function in $n+m$ variables with $n,m \in \N$. Let $\vec{\zeta} \in \C^{+m}$ be fixed and let $\til{q}_{\vec{\zeta}}$ be the the \HN function of $n$ variables defined as
$$\til{q}_{\vec{\zeta}}\colon \vec{z} \mapsto q(\vec{z},\vec{\zeta}).$$
Then, the following statements hold.
\begin{itemize}
    \item[(i)]{The function
    $$\vec{\zeta} \mapsto \int_{\R^n}\psi(\vec{t})\diff\til{\mu}_{\vec{\zeta}}(\vec{t})$$
    is pluriharmonic on $\C^{+m}$ for any function $\psi$ as in the Stieltjes inversion formula \eqref{eq:Stieltjes_inversion}.
    }
    \item[(ii)]{The function
    $$\vec{\zeta} \mapsto \til{\mu}_{\vec{\zeta}}(U)$$
    is pluriharmonic on $\C^{+m}$ for any Borel set $U \subseteq \R^n$. In particular, let $(\til{\mu}_{\vec{\zeta}})_\mathrm{a.c.}$ and $(\til{\mu}_{\vec{\zeta}})_\mathrm{sing}$ be the absolutely continuous and singular part of $\til{\mu}_{\vec{\zeta}}$ with respect to $\lambda_{\R^n}$ in accordance with the Lebesgue decomposition theorem. Then, the functions
    $$\vec{\zeta} \mapsto (\til{\mu}_{\vec{\zeta}})_\mathrm{a.c.}(U) \quad\text{and}\quad \vec{\zeta} \mapsto (\til{\mu}_{\vec{\zeta}})_\mathrm{sing}(U)$$
    are pluriharmonic on $\C^{+m}$ for any Borel set $U \subseteq \R^n$.
    }
\end{itemize}
\end{prop}

\proof 
This is obtained as a consequence of  Theorem \ref{thm:singular_part_polydisc}(ii) and Lemma \ref{lem:measure_polydisc} employing  the connection between the representing measures for RP-functions on the polydisc and Nevanlinna functions on the poly-upper half-plane described in Section \ref{subsec:correspondence}.
\endproof

\section{Estimates on the measure of cubes}
\label{sec:estimates}

We begin by recalling from \cite[Sec. 13.3]{Vladimirov1979} the universal upper bound for a \HN function. For a function $\til{q}$ of one variable, this says that there exists a constant $M$ such that the estimate
\begin{equation}
    \label{eq:HN_estimate_1var}
    |\til{q}(z)| \leq M\,\frac{1+|z|^2}{y}
\end{equation}
holds for all $z = x + \I\,y \in \C^+$. If we instead have $q$ which is a \HN function of several variables, the above estimate takes the form
\begin{equation}
    \label{eq:HN_estimate_Nvar}
    |q(\vec{z})| \leq M\,\frac{1+\sum_{j=1}^n|z_j|^2}{\left(\sum_{j=1}^ny_j^2\right)^{\frac{1}{2}}},
\end{equation}
where $\vec{z} = \vec{x} + \I\,\vec{y}\in \C^{+n}$.

Furthermore, denote by
$$D(\vec{\tau},y) := \big\{\vec{t} \in \R^n~\big|~\forall\,j=1\ldots,n:|\tau_j - t_j| < y\big\}$$
the open polydisc in $\R^n$ with centre $\vec{\tau}$ and uniform radius $y$ . We have the following estimate:

\begin{prop}
\label{prop:volume_estimate}
Let $\mu$ be a Nevanlinna measure. Then, there exists a constant $M \geq 0$ such that for any $\vec{\tau} \in \R^n$ and $y > 0$, we have
\begin{equation}
    \label{eq:volume_estimate}
    \mu(D(\vec{\tau},y)) \leq \frac{M}{(2\pi)^n}\,\frac{y^{n-1}}{\sqrt{n}}(1+\|\vec{\tau}\|_2^2+ny^2).
\end{equation}
\end{prop}

\proof
Let $q$ be the \HN function given by $(0,\vec{0},\mu)$ in the sense of Theorem \ref{thm:intRep_Nvar} and let $v$ be its imaginary part. Using \cite[Cor. 4.6(ii)]{LugerNedic2019}, we estimate that
\begin{multline*}
    v(\vec{\tau} + \I\,y\,\vec{1}) = \frac{1}{\pi^n}\int_{\R^n}\pois_n(\vec{\tau} + \I\,y\,\vec{1},\vec{t})\diff\mu(\vec{t}) = \frac{1}{\pi^n}\int_{\R^n}\prod_{j=1}^n\frac{y}{(\tau_j - t_j)^2+y^2}\diff\mu(\vec{t}) \\ \geq \frac{1}{\pi^n}\int_{D(\vec{\tau},y)}\prod_{j=1}^n\frac{y}{(\tau_j - t_j)^2+y^2}\diff\mu(\vec{t}) \geq \frac{1}{\pi^n}\int_{D(\vec{\tau},y)}\prod_{j=1}^n\frac{1}{2y}\diff\mu(\vec{t}) \\ \geq \frac{1}{(2\pi y)^n}\mu(D(\vec{\tau},y)).
\end{multline*}
On the other hand, using estimate \eqref{eq:HN_estimate_Nvar}, we deduce that
$$v(\vec{\tau} + \I\,y\,\vec{1}) \leq |q(\vec{\tau} + \I\,y\,\vec{1})| \leq M\,\frac{1+\sum_{j=1}^n|\tau_j + \I\,y|^2}{\left(\sum_{j=1}^ny^2\right)^{\frac{1}{2}}} = M\,\frac{1 + \|\vec{\tau}\|_2^2 + ny^2}{\sqrt{n}y}.$$
Combining these two estimates yields  the desired result.
\endproof

Note that if
$$B(\vec{\tau},y) := \big\{\vec{t} \in \R^n~\big|~\|\vec{\tau}- \vec{t}\|_2^2 < y^2\big\}$$
is the open ball in $\R^n$ with centre $\vec{\tau}$ and radius $y$, then $\mu(B(\vec{\tau},y)) \leq \mu(D(\vec{\tau},y))$. Hence, the right-hand side of inequality \eqref{eq:volume_estimate} also provides an upper estimate for the measure of a ball in $\R^n$.

Consider now the following example illustrating  the estimate \eqref{eq:volume_estimate} for different measures in $\R^2$.

\begin{example}
The Dirac measure $\delta_{(0,0)}$ at $\vec{0} \in \R^2$ is known not be a Nevanlinna measure and clearly does  not satisfy estimate \eqref{eq:volume_estimate}.  On the other hand,  the positive Borel measure on $\R^2$ 
$$\mu_2(U) = \pi\int_\R\chi_U(t,-t)\diff t$$
is a Nevanlinna measure, see \eg \cite[Ex. 3.14]{LugerNedic2021} or Corollary \ref{coro:measure_on_lines} and,  satisfies estimate \eqref{eq:volume_estimate}. One may note that the Hausdorff measure on the main diagonal satisfies the estimate, but as we already have noted, fails to be a Nevanlinna measure.\hfill$\lozenge$
\end{example}

Observe now that the  estimate \eqref{eq:volume_estimate} implies for large $y$
$$0 \leq \frac{\mu(D(\vec{\tau},y))}{y^{n+1}} \leq M_\infty.$$
In other words, the measure of a cube cannot grow faster than some constant times $y^{n+1}$. On the other hand, if there happened to exists $y_1 < y_2$ such that $0 < \mu(D(\vec{\tau},y_1)) = \mu(D(\vec{\tau},y_2))$, then $\mu$ would be a non-trivial finite Nevanlinna measure, contradicting Corollary \ref{coro:finite_measure}. Hence, it is natural to inquire for the slowest rate of growth, though at this point we cannot be certain whether the growth rate could arbitrarily slow. Consider, hence, the following proposition. 

\begin{prop}
\label{prop:limsup}
Let  and let $\mu$ be a non-trivial Nevanlinna measure on $\R^2$. Then, it holds that 
\begin{equation}
    \label{eq:lower_bound}
  \limsup_{R\to\infty} R^{-1}\mu(D(\vec{0},R))>0.
\end{equation}
\end{prop}

\proof
Let $q$ be the \HN function on $\C^{+2}$ given by the data $(0,\vec{0},\mu)$ in the sense of Theorem \ref{thm:intRep_Nvar}. Denote $\widetilde{q}(z) : =q(z,z)$ for $z\in \C^+$, so that $\widetilde{q}$ is a \HN function in one variable. Since
$$\Im[\widetilde{q}(\I)] = \Im[q(\I,\I)] = \frac{1}{\pi^2}\int_{\R^2}\frac{1}{1+t_1^2}\,\frac{1}{1+t_2^2}\diff\mu(\vec{t}) >0$$
due to $\mu$ being non-trivial, $\Im[\widetilde{q}]$ does not vanish identically on $\C^+$. Hence, it holds that
\begin{equation}\label{eq:limit}
\lim_{y\to\infty}y\,\Im[\widetilde{q}(\I\,y)] > 0,
\end{equation}
where the limit may also take the value $\infty$. Namely, if $\widetilde{q}$ is represented by $(\widetilde{a},\widetilde{b},\widetilde{\mu})$, the above limit equals $+\infty$ if $\widetilde{b} > 0$, while in case $\widetilde{b} = 0$, Lebesgue's monotone convergence theorem verifies that the limit equals $\mu(\R)$ (whether finite or infinite) \cite{KacKrein1974}. On the other hand, using representation \eqref{eq:intRep_Nvar} to describe $\Im[q]$, we derive that
\begin{equation*}
y\,\Im[\widetilde{q}(\I\,y)]  = \frac{1}{\pi^2}\,\frac{1}{y}\int_{\R^2}\frac{y^2}{t_1^2+y^2}\,\frac{y^2}{t_2^2+y^2}\diff\mu(\vec{t}).
\end{equation*}

Take now $y=2^{\ell_0}$, where $\ell_0\in \N$, and denote $D_k := D(\vec{0},2^k)$. Divide the domain of integration in the above integral to the cube $D_{\ell_0}$ and the "annuli" $D_k-D_{k-1}$, for $k>\ell_0$ to get
\begin{multline}
\label{eq:estimate1}
    2^{\ell_0}\,\Im[\widetilde{q}(2^{\ell_0})] = \pi^{-2}\,2^{-\ell_0}
    \cdot \left(\int_{D_{\ell_ 0}}\frac{2^{2\ell_0}}{t_1^2+2^{2\ell_0}}\,\frac{2^{2\ell_0}}{t_2^2+2^{2\ell_0}}\diff\mu(\vec{t}) \right. \\ \left. + \sum_{k=\ell_0 + 1}^\infty\int_{D_{k}-D_{k-1}}\frac{2^{2\ell_0}}{t_1^2+2^{2\ell_0}}\,\frac{2^{2\ell_0}}{t_2^2+2^{2\ell_0}}\diff\mu(\vec{t})\right).
\end{multline}
In the first integral, we note that $0 \leq |t_j| \leq 2^{\ell_0}$ for both $j = 1,2$ and we hence bound the integrand form above by $1$, yielding the estimate
\begin{equation}
\label{eq:cube_estimate1}
\int_{D_{\ell_ 0}}\frac{2^{2\ell_0}}{t_1^2+2^{2\ell_0}}\,\frac{2^{2\ell_0}}{t_2^2+2^{2\ell_0}}\diff\mu(\vec{t}) \leq \mu(D_{\ell_0}).
\end{equation}
In the second integral, there is always at least one index $j = 1,2$ for which $2^{k-1} < |t_j| \leq 2^k$. For one such index $j$, we make the estimate that
$$\frac{2^{2\ell_0}}{t_j^2+2^{2\ell_0}} \leq \frac{2^{2\ell_0}}{2^{2(k-1)}+2^{2\ell_0}} = \frac{1}{2^{2(k-\ell_0-1)}+1} \leq 2^{-2(k-\ell_0-1)},$$
while for all other index $j'$ we use the estimate
$\frac{2^{2\ell_0}}{t_{j'}^2+2^{2\ell_0}} \leq 1.$
Hence
\begin{multline}
    \label{eq:cube_estimate2}
    \int_{D_{k}-D_{k-1}}\frac{2^{2\ell_0}}{t_1^2+2^{2\ell_0}}\,\frac{2^{2\ell_0}}{t_2^2+2^{2\ell_0}}\diff\mu(\vec{t}) \\ \leq 2^{-2(k-\ell_0-1)}\mu(D_{k}\setminus D_{k-1}) \leq 2^{-2(k-\ell_0-1)}\mu(D_{k}).
\end{multline}
Using estimates \eqref{eq:cube_estimate1} and \eqref{eq:cube_estimate2} to obtain an upper bound for the right-hand side of equality \eqref{eq:estimate1} yields
$$\begin{array}{RCL}
    \multicolumn{3}{L}{2^{\ell_0}\,\Im[\widetilde{q}(2^{\ell_0}\,\I)]} \\
    ~ & \leq & \pi^{-2}\,2^{-\ell_0}\left(\mu(D_{\ell_0}) +  \sum_{k=\ell_0+1}^\infty 2^{-2(k-\ell_0-1)}\mu(D_{k})\right) \\[0.4cm]
    ~ & = & \pi^{-2}\left(2^{-\ell_0}\mu(D_{\ell_0}) + 4\sum_{k=\ell_0+1}^\infty 2^{-(k-\ell_0)}\,2^{-k}\mu(D_{k})\right) \\[0.4cm]
    ~ & \leq & 4\,\pi^{-2} \sum_{k=\ell_0}^\infty 2^{-(k-\ell_0)}\,2^{-k}\mu(D_{k}) \leq 4\,\pi^{-2}\,\sup_{k\geq \ell_0}\big(2^{-k}\mu(D_{k})\big).
\end{array}$$
The above estimate shows that if estimate \eqref{eq:limit} would be false, \ie if 
$$\limsup_{k\to\infty}\big(2^{-k}\mu(D_k)\big) = 0,$$
then it would hold that $\lim_{\ell\to\infty}2^\ell\,\Im[\widetilde{q}(2^\ell\,\I)] =0$, contradicting \eqref{eq:limit}.
\endproof

\begin{remark} The above proof can easily be adapted for Nevanlinna measures in dimensions $n \geq 3$, yielding the same estimate. However, based on various examples, \eg $q(z)=-(z_1+\ldots+z_n)^{-1}$ in which the Nevanlinna measure is the $(n-1)$-dimensional Hausdorff measure restricted to the hyper-plane $\{\vec{t} \in \R^n~|~t_1+\ldots +t_n = 0\}$ \cite[Thm. 4.2]{Nedic2019}, we conjecture that the true order of growth in dimension $n$ is 
$$
  \limsup_{R\to\infty} R^{1-n}\mu(D(\vec{0},R))>0.
$$
One may also ask whether above the limit always exists with value in $(0,\infty]$.
\end{remark}

\bibliographystyle{amsplain}

\end{document}